\documentclass [10pt,reqno]{amsart}
\usepackage{ucs}
\usepackage{amsmath, amssymb, amscd}
\usepackage{pb-diagram}
\usepackage[latin1]{inputenc}
\usepackage{amsmath}
\usepackage{amssymb}

\newtheorem{theorem}{Theorem}[section]
\newtheorem{definition}[theorem]{Definition}
\newtheorem{lemma}[theorem]{Lemma}
\newtheorem{prop}[theorem]{Proposition}
\newtheorem{corollary}[theorem]{Corollary}

\newtheorem*{remark}{Remark}

\newcommand{\abs}[1]{\lvert#1\rvert}
\newcommand{\norm}[1]{\| #1 \|}
\newcommand{\klammer}[1]{\left( #1 \right)}
\newcommand{\eckklammer}[1]{\left[ #1 \right]}

\newcommand{\geklatau}[1]{\left\{ #1 \right\}_\tau}
\newcommand{\dbar}{\bar{\partial}}
\renewcommand{\epsilon}{\varepsilon}

\DeclareMathAlphabet{\mathpzc}{OT1}{pzc}{m}{it}
\usepackage{mathrsfs}

\newcommand{\dirac}{\diagup \hskip-9pt D}

\newcommand{\Z}{\mathbb{Z}}
\newcommand{\C}{\mathbb{C}}

\newcommand{\R}{\mathbb{R}}

\renewcommand{\qed}{$\hfill \square$ \bigskip \\}
\renewcommand{\phi}{\varphi}

\newcommand{\ad}{\text{ad}}

\newcommand{\rk}{\text{rk}}

\newcommand{\tensor}{\otimes}
\newcommand{\End}{\text{End}}
\newcommand{\Hom}{\text{Hom}}
\newcommand{\Aut}{\text{Aut}}
\renewcommand{\det}{\text{det}}
\newcommand{\id}{\text{id}}

\newcommand{\ind}{\text{ind}}

\newcommand{\Atheta}{\mathscr{A}_{\theta}(P_E)}

\newcommand{\su}{\mathfrak{su}}
\renewcommand{\sl}{\mathfrak{sl}}
\newcommand{\gl}{\mathfrak{gl}}
\renewcommand{\u}{\mathfrak{u}}

\renewcommand{\bar}{\overline}
\newcommand{\G}{\mathscr{G}^0}

\newcommand{\tr}{\text{tr}}
\newcommand{\data}{\mathfrak{s},E}
\newcommand{\dataf}{\mathfrak{s},F}
\newcommand{\conf}{\mathscr{C}}
\newcommand{\bonf}{\mathscr{B}}
\newcommand{\s}{\mathfrak{s}}

\begin{document}

\thispagestyle{empty}

\title[On Higher rank instantons \& the Monopole cobordism program]{On Higher rank instantons \& the Monopole cobordism program}
 \author
[Raphael
Zentner]{Raphael Zentner} 

\begin {abstract} 
Witten's conjecture suggests that the polynomial invariants of Donaldson are expressible in terms of the Seiberg-Witten invariants if the underlying four-manifold is of simple type.  
A higher rank version of the Donaldson invariants was introduced by Kronheimer. Before even having been defined, the physicists Mari\~no and Moore had already suggested that there should be a generalisation of Witten's conjecture to this type of invariants. We study a generalisation of the classical cobordism program to the higher rank situation and obtain vanishing results which gives evidence that the generalisation of Witten's conjecture should hold.
\end {abstract}

\address {Mathematisches Institut \\ Westf\"alische Wilhelms-Universit\"at \\ Einsteinstr. 62 \\ 48149 M\"unster\\
Germany}

\email{raphael.zentner@math.uni-muenster.de}

\maketitle

\section*{Introduction}
Recently Peter Kronheimer introduced polynomial invariants associated to moduli spaces of anti-selfdual $PU(N)$ connections in Hermitian rank-N-bundles \cite{K}. These are generalisations of the Donaldson invariants, but the technical problems are much harder than in the classical situation. Before these invariants were even properly defined, the physicists Mari\~no and Moore \cite{MM} had conjectured that there should be a generalisation of Witten's conjecture for these invariants, implying in particular that they do not contain new differential-topological information. Their argument relies again on physics. Kronheimer computed explicitely his higher rank invariants for manifolds obtained through knot surgery on the K3-surface, confirming the conjecture for this class of examples.
\\

The main important gauge-theoretical invariants of a smooth closed four-manifold are the polynomial invariants of Donaldson \cite{D}, derived from anti-selfdual $PU(2)$ connections in rank-2-bundles, and the Seiberg-Witten invariants, derived from the Seiberg-Witten equations \cite{W} which are associated to $Spin^c$ structures on the four-manifold. Kronheimer and Mrowka have proved a structure theorem for the Donaldson invariants \cite{KM2}, showing that for 4-manifolds of `simple type' the polynomial invariants are specified by certain algebraic-topological data, in particular the intersection form, and a finite set of distinguished cohomology classes in the group $H^2(X,\Z)$ each coming with some rational coefficient. With this at hand, Witten claims that the polynomial invariants are determined by the Seiberg-Witten-invariants, with the basic classes being the first Chern classes of the $Spin^c$-structures with non-trivial Seiberg-Witten-invariant, and with an explicit formula for
the rational coefficients \cite{W}. 

Witten derived this conjecture from correlation functions in quantum field theory and certain limiting behaviours with respect to a certain coupling constant. Mathematicians then tried to derive a proof of the conjecture by a certain cobordism obtained from $PU(2)$ monopoles. Heuristically the idea is as follows. There is a circle action on the moduli space of $PU(2)$ monopoles which comes from multiplying the spinor component by a complex number of norm one. The fixed point set of this action consists of the moduli space of $PU(2)$ instantons, and further a finite number of moduli spaces of classical $U(1)$ Seiberg-Witten monopoles. The circle acts freely on the complement of this fixed point locus, and the quotient yields a cobordism between a projective bundle over the moduli space of $PU(2)$ instantons, and projective bundles over the moduli spaces of $U(1)$ Seiberg-Witten monopoles. 

Furthermore, the canoncial cohomology class that yields the polynomial invariant of Donaldson extends to the cobordism. The evaluation of this extension on one side, yielding the Donaldson invariant, is therefore equal to the evaluation on the other sides, which should be expressions containing the Seiberg-Witten invariants of the moduli spaces in the fixed point locus. This program was started independently by Pidstrigach and Tyurin \cite{PT}, and Okonek and Teleman \cite{OT2}, \cite{T}. It was carried on over years by Feehan and Leness \cite{FL1,FL2,FL3,FL4}. It seems that they have now proved the full conjecture \cite{FL5}.
\\

Our intention is to consider the generalisation of the cobordism program with the perspective of applying it to the generalisation of Witten's conjecture by  Mari\~no and Moore. We introduce $PU(N)$ monopoles: For a given $Spin^c$ structure $\mathfrak{s}$ and a Hermitian rank-N-bundle $E$ on $X$, the configuration space will consist of sections $\Psi$ of the `twisted spinor bundle' $W^+ = S^+_\mathfrak{s} \tensor E$ and by unitary connections $A$ in $E$ with fixed induced connection $\theta$ in the determinant line bundle of $E$. The straightforward generalisation of the $PU(2)$ monopole equations then read:
\begin{equation*}
\begin{split}
 \dirac_{A}^{+} \Psi & =  0 \\
 \gamma((F_A^+)_{0}) - \mu_{0,0}(\Psi) & =  0 \ .
\end{split}
\end{equation*}
Here $\dirac_{A}^{+}$ is the associated Dirac-operator to $A$, the map $\gamma$ is derived from Clifford-multiplication, $(F_A^+)_{0}$ is the self-dual part of the curvature of the $PU(N)$ connection $A$ induced by $A$, and $\mu_{0,0}$ is a quadratic map in the spinor which is explicitely described below. The gauge group of the problem is that of special unitary automorphisms of $E$. 

Again, there is a circle action on the moduli space of these $PU(N)$ monopoles which is given by the formula $(z,[\Psi,A]) \mapsto [z^{1/N} \Psi,A]$. The moduli space of $PU(N)$ instantons is contained as the locus of monopoles with vanishing spinor, and the other fixed point loci are labelled by a finite number of isomorphism classes of proper subbundles $[F]$ of $E$. An equivalence class $[\Psi,A]$ belongs to the $[F]$-locus $M^{[F]}$ if for each $F \in [F]$ there is a representative $(\Psi,A)$ with the spinor $\Psi$ being a section of $S^+_\mathfrak{s} \tensor F$ and with the connection $A$ keeps the proper subbundle $F$ invariant. It turns out that if $X$ is simply connected the description of $M^{[F]}$ is particularly simple after fixing one such $F$: The content of Theorem \ref{S1fixedpointset_main_text} is that we have a `parametrisation' 
\begin{equation}\label{contributions}
 M_{\mathfrak{s},F} \times M^{asd}_{F^\perp} \to M^{[F]} \ .
\end{equation}
Here $M_{\mathfrak{s},F}$ is a moduli space of $U(n)$ monopoles with $n= \rk(F)$ having possible values $1 \leq n < N$, and $M^{asd}_{F^\perp}$ is the moduli space of anti-selfdual $PU(N-n)$ connections in $F^\perp$. In the case $n=1$ the moduli space $M^U_{\mathfrak{s},F}$ is a classical $U(1)$ Seiberg-Witten moduli space. The map (\ref{contributions}) is surjective and is bijective if restricted to the subspaces of the corresponding moduli spaces which consist of elements with zero-dimensional stabiliser.
In the classical case $N=2$ we can only have $n=1$ and there are no non-trivial $PU(1)$ connections. 

Now the components $M^{[F]}$ are the possible contributions to the formula expressing the $PU(N)$ instanton invariant according to the cobordism program indicated above. The generalisation of Witten's conjecture to the $PU(N)$ instanton invariants would follow if only those components $M^{[F]}$ contribute in a non-trivial way for which we have $n=\rk(F)=1$. We give two results that indicate that this should be true - one general but under slightly speculative simplifying assumption, the other on K\"ahler surfaces. 
\\

The first section sets up our configuration space, introduces the above mentioned quadratic map $\mu_{0,0}$ and variations $\mu_{0,\tau}$ of it with a parameter $\tau \in [0,1]$. We derive some important properness property of it. The $PU(N)$ monopole equations are then introduced and it is indicated how to obtain an Uhlenbeck-type compactification of the moduli space. With only minor modifications, we then introduce the $U(n)$ monopole equations. The second section studies the circle action on the moduli space of $PU(N)$ monopoles, analyses its fixed point set and relates it to $U(n)$ monopoles and $PU(N-n)$ instantons. The third section gives our first argument why we expect no contributions from fixed point loci $M_{\s,F} \times M^{asd}_{F^{\perp}}$ for $\rk(F) > 1$. We prove that under certain assumptions the moduli space $M_{\s,F}$ is cobordant to the empty space. For the second argument, we analyse the $U(n)$ monopole equations on K\"ahler surfaces in the forth section. By studying decoupling phenomena, we show that the moduli space becomes empty as soon as we perturb by a non-zero holomorphic two-form. 

\section*{Acknowledgements} This manuscript constitutes a part of the author's PhD thesis. He is grateful to his former advisor Andrei Teleman for the care taken in guiding him and the many mathematical discussions involved 
. He is also indebted to Peter Kronheimer and Kim Fr\o yshov for discussions on some of the aspects.

\section{Preliminaries, $PU(N)$ and $U(n)$ monopoles}
Here we shall introduce the $PU(N)$ monopole as well as the $U(n)$ monopole equations associated to the data
of a $Spin^c$-structure $\mathfrak{s}$ and a Hermitian bundle $E \to X$ on a Riemannian four-manifold $X$.
 We shall define the monopole equations, the moduli space, and prove a uniform bound on the spinor
component of a solution to the monopole equations.  We assume knowledge with standard (abelian) Seiberg-Witten theory as for instance found in the monographs \cite{N,M} or the lecture notes \cite{T3}.

\subsection{Twisted $Spin^{c}$ structures and associated Dirac operators}
Let $X$ be a closed oriented Riemannian four-manifold with a
$Spin^c$ structure $\mathfrak{s}$ on it. The $Spin^c$ structure consists of two
Hermitian rank 2 vector bundles $S^\pm_{\mathfrak{s}}$ with identified
determinant line bundles and a Clifford multiplication 
\begin{equation*} 
\gamma : \Lambda^1(T^*X) \to \Hom_\C(S^+_\mathfrak{s},S^-_\mathfrak{s}) \ .
\end{equation*}

The Clifford map $\gamma$ is, up to a universal constant, an isometry of the
cotangent bundle
onto a real form inside $\Hom_\C(S^+_\mathfrak{s},S^-_\mathfrak{s})$ which can
be specified by the Pauli matrices. We extend $\gamma$ to $\End(S^+_\s
\oplus S^-_\s)$ by $ - \gamma^*$ on the negative Spinor bundle. It then
naturally extends to exteriour powers of $T^*X$, and in particular its
restriction to self-dual two-forms is zero on the negative Spinor bundle, and
induces an isomorphism 
\begin{equation*}
\gamma: \Lambda^2_+(T^*X) \stackrel{\cong}{\to} \su(S^+_\mathfrak{s}) \ .
\end{equation*}
Furthermore suppose we are given a
Hermitian vector bundle $E$ with determinant line bundle $w = \det(E)$ on $X$. We can then form spinor bundles 
\begin{equation*}
W_{\data}^\pm := S^\pm_{\mathfrak{s}} \tensor E .
\end{equation*}
Clifford multiplication extends by tensoring with the identity on $E$. 
This way we obtain a $Spin^c$ - structure `twisted' by the hermitian bundle $E$.
\\

Taking the tensor product with the identity on $E$ induces a Clifford
multiplication $\gamma : \Lambda^1(T^*X) \to \Hom(W^+_{\data},W^-_{\data})$.
Let's fix a background $Spin^c$ connection $B$ on $\mathfrak{s}$ (that will always stay the same) and suppose we are given a unitary connection $A$ on $E$. By composing
the tensor product connection $\nabla_{B} \tensor \nabla_{A}$ with the Clifford
multiplication we get a Dirac operator
\begin{equation*}
\dirac_{A}^{\pm} := \gamma \circ \left(\nabla_{B} \tensor \nabla_{{A}}\right) :
\Gamma(X;W^\pm_{\data}) \to \Gamma(X;W^\mp_{\data}) \ .
\end{equation*}
This is a self-adjoint first order elliptic operator. We have
oppressed the $Spin^c$ connection $B$ from the notation because it will not be a variable in our theory.


\subsection{Algebraic preliminaries} We shall now introduce the quadratic map $\mu_{0,\tau}: S^+\tensor E \to \su(S^+)\tensor_\R \su(E)$, defined for a real number $\tau \in [0,1]$, that appears in the non-abelian monopole equations. For $\tau = 0$ it will
appear in the $PU(n)$ - monopole equations, and for $\tau \neq 0$ in the $U(n)$ monopole equations. This map is a natural generalisation of the corresponding map appearing in the situation of $PU(2)$ monopoles \cite{FL2, FL4, T2} as well as the one in the classical (abelian) Seiberg-Witten equations \cite{KM,W}. 

The twisted spinor bundles $W^\pm_{\data}$ are associated bundles of the fibre product of a
$Spin^c$
principal bundle and a $U(n)$-principal bundle on X, with the standard fibre
$\C^2 \tensor \C^n$.
\\

Let us consider the isomorphism
\begin{equation*}
\begin{split}
  (p,q): \ \,  \gl(\C^n) & \to \sl(\C^n) \oplus \C \, \id \\
   	 a  & \mapsto  \left(a - \frac{1}{n} \ \tr(a) \cdot \id, \frac{1}{n} \
\tr(a) \cdot \id\right) \ .
\end{split}
\end{equation*}
Both components $p$ and $q$ are orthogonal projections onto their images. Note
that $\gl(\C^2) \tensor \gl(\C^n)$ and $\gl(\C^2 \tensor \C^n)$ are canonically
isomorphic. We define the orthogonal projections
\begin{equation*}
\begin{split}
  P: \ \gl(\C^2 \tensor \C^n) & \to \sl(\C^2) \tensor \sl(\C^n) \ , \\
  Q: \ \gl(\C^2 \tensor \C^n) & \to \sl(\C^2) \tensor \C \, \id
\end{split}
\end{equation*}
to be the tensor product $( \ )_0 \tensor p$ respectively $( \ )_0 \tensor q$,
with $( \ )_0$ denoting the trace-free part of the endomorphism of the first
factor $\C^2$. 
 
For elements $\Psi, \Phi \in \C^2 \tensor \C^n$ we define
\[
\mu_{0,\tau}(\Psi,\Phi):=P (\Psi \Phi^*) \ + \tau \, Q ( \Psi \Phi^*) ,
\]
 where $(\Psi \Phi^*) \in \gl(\C^2 \tensor \C^n)$ is defined to be the
endomorphism $\Xi \mapsto
\Psi (\Phi,\Xi)$. 

With this notation $\mu_{0,1} (\Psi,\Phi)$ is simply the orthogonal projection
of the endomorphism $\Psi \Phi^* \in \gl(\C^2 \tensor \C^n)$ onto $\sl(\C^2)
\tensor \gl(\C^n)$. We shall also write $\mu_{0,\tau}(\Psi):=
\mu_{0,\tau}(\Psi,\Psi)$ for the associated quadratic map. In the case $n=1$ the
map $\mu_{0,1}(\Psi)$ is the quadratic map in the spinor usually occuring in the
Seiberg-Witten equations \cite{W} \cite{KM}. The proof of the following proposition can be found in \cite{Z_un}. \\
%

\begin{prop}\label{properness}
Suppose $n > 1$. Then the quadratic map $\mu_{0,\tau}$ is uniformly proper. In
other words, there is a positive
constant $c > 0$ such that 
  \begin{equation}\label{propernessconstant}
    \abs{\mu_{0,\tau}(\Psi)} \geq c \abs{\Psi}^2 \ .
  \end{equation}
As a consequence we have the formula
  \begin{equation}
    \left(\mu_{0,\tau}(\Psi) \Psi, \Psi \right) \geq c^2 \abs{\Psi}^4 \ 
  \end{equation}
whenever $\tau \geq 0$.
Furthermore, suppose $n \geq 2$ or $\tau \neq 0$. 
Then the bilinear map $\mu_{0,\tau}$ is `without zero-divisors' in the following sense: If
$\mu_{0,\tau}(\Psi,\Phi) = 0$, then either $\Psi = 0$ or $\Phi = 0$.
\end{prop}

{\em Proof:}
Obviously we have
$\abs{\mu_{0,\tau}(\Psi)} \geq \abs{\mu_{0,0}(\Psi)}$, so for the first
assertion it will be enough to consider $\mu_{0,0}$ alone. We will show that
$\mu_{0,0}(\Psi) =
0$ implies $\Psi = 0$. Because $\mu_{0,0}$ is quadratic and the unit sphere
inside
$\C^2 \tensor \C^n$ is compact we then get the claimed uniform
properness-inequality (\ref{propernessconstant}). 

We shall use the canonical
isomorphism $\C^2
\tensor \C^n \cong \C^n \oplus \C^n$, which permits to write a general element 
\[
  \Psi=\begin{pmatrix} 1 \cr 0 \end{pmatrix} \tensor \alpha + 
      \begin{pmatrix} 0 \cr 1 \end{pmatrix} \tensor \beta 
\]
as 
\[
  \Psi= \begin{pmatrix}\alpha \cr \beta \end{pmatrix} \ .
\]
We then have
\begin{equation*}
\begin{split}
\mu_{0,0}(\Psi) & =  P (\Psi \Psi^*) \\
	& =  P \left( \begin{pmatrix} \alpha \cr \beta \end{pmatrix}
\begin{pmatrix} \alpha^* & \beta^*\end{pmatrix} \right) \\
	& =  P \begin{pmatrix} \alpha \alpha^* & \alpha \beta^* \cr 
			\beta \alpha^* & \beta \beta^* \end{pmatrix} \\
	& =   \begin{pmatrix} \frac{1}{2} (\alpha \alpha^* - \beta \beta^*)_0 & 
	      (\alpha \beta^*)_0 \cr (\beta \alpha^*)_0 & \frac{1}{2}(\beta
\beta^* - \alpha \alpha^*)_0 \end{pmatrix} \ . 
\end{split}
\end{equation*}
In particular, if $\mu_{0,0}(\Psi)= 0$, then we have $(\alpha \beta^*)_0 = 0$. 
\begin{lemma} \label{lemmainlemma}
The equation $(\alpha \beta^*)_0 = 0$ implies that $\alpha=0$ or $\beta=0$. In
other words, the bilinear map $(\alpha, \beta) \to (\alpha \beta^*)_0$ is
without
zero-divisors (here $n \geq 2$ is implicitely understood).
\end{lemma}
{\em Proof of Lemma \ref{lemmainlemma}:} Write the elements $\alpha$ and $\beta
$ as 
\[
  \alpha = (\alpha_i)_{i=1}^n \ , \quad \beta= (\beta_i)_{i=1}^n \ .
\]
Then the equation $(\alpha \beta^*)_0 = 0 $ reads in  matrix-form
\begin{equation*}
\begin{split}
\left( \begin{array}{ccc} \alpha_1 \bar{\beta_1} - \frac{1}{n} \sum \alpha_i
\bar{\beta_i} & \dots & \alpha_1 \bar{\beta_n} \\
\vdots & \ddots & \vdots \\
\alpha_n \bar{\beta_1} & \dots & \alpha_n \bar{\beta_n}-\frac{1}{n} \sum
\alpha_i \bar{\beta_i} \end{array} \right) = 0  \ .
\end{split}
\end{equation*}
Suppose $\beta\neq 0$, for instance $\beta_j \neq 0$. Then the $j^{\text{th}}$
column implies that $\alpha_i= 0$ for all $i \neq j$. Thus the
$j^{\text{th}}$ element in the $j^{\text{th}}$ column simplifies,
\begin{equation*}
 \alpha_j \bar{\beta_j} - \frac{1}{n}\sum_{i=1}^{n} \alpha_i \bar{\beta_i} \, =
\, (1-\frac{1}{n}) \alpha_j \bar{\beta_j} \ .
\end{equation*}
Therefore we have $\alpha_j= 0$ as well, so that we have $\alpha=0$. The case
$\alpha \neq 0 $
is analogous. \qed
Returning to the problem
\[
  \begin{pmatrix} \frac{1}{2} (\alpha \alpha^* - \beta \beta^*)_0 & 
	      (\alpha \beta^*)_0 \cr (\beta \alpha^*)_0 & \frac{1}{2}(\beta
\beta^* - \alpha \alpha^*)_0 \end{pmatrix} = 0 \ ,
\]
we see that the lemma gives $\alpha = 0$ or $\beta = 0$. Suppose, without loss
of
generality, that the first is the case. Then we are left with $(\beta \beta^*)_0
= 0 $. Now again with lemma \ref{lemmainlemma} we see that this also implies
$\beta=
0$. Therefore 
\[
\Psi= \begin{pmatrix} \alpha \cr \beta \end{pmatrix} = 0 \ .
\]

The second assertion now follows from the first, remembering that $P$ and $Q$
are both orthogonal projections. For non-negative $\tau$ we have the
inequality 
\begin{equation*}
\begin{split}
\left(\mu_{0,\tau}(\Psi)\Psi,\Psi \right) = & \, \left(P (\Psi \Psi^*) \Psi,
\Psi
\right) + \tau \klammer{Q(\Psi \Psi^*) \Psi, \Psi} \\
  = & \, \left( P (\Psi \Psi^*), \Psi \Psi^* \right) + \tau \klammer{Q(\Psi 
\Psi^*), \Psi \Psi^*} \\
  = & \, \left( P (\Psi \Psi^*), P (\Psi \Psi^*) \right) + \tau \klammer{Q(\Psi
\Psi^*), Q(\Psi \Psi^*)} \\
  \geq & \ \abs{\mu_{0,0}(\Psi)}^2 \\
  \geq & \ c^2 \abs{\Psi}^4 \ .
\end{split}
\end{equation*}

Similarly, one shows the claimed property about the `zero-divisors'. 
\qed 

Because of the equivariance property of the map $\mu_{0,\tau}$ we get in
a straightforward way corresponding maps between bundles, giving rise to
\begin{equation*}
\mu_{0,\tau}: W^\pm_{\data}  \times  W^\pm_{\data} \ 
\to \sl(S^\pm_{\s}) \tensor_{\C} \gl(E) \ ,
\end{equation*}
respectively, for the quadratic map,
\begin{equation*}
\mu_{0,\tau}: W^\pm_{\data} \ 
\to \su(S^\pm_{\s}) \tensor_{\R} \u(E) \ , 
\end{equation*}
or
\begin{equation*}
\mu_{0,0}: W^\pm_{\data} \ 
\to \su(S^\pm_{\s}) \tensor_{\R} \su(E) \ 
\end{equation*}
if we put $\tau = 0$. 

These maps on the bundle level satisfy the corresponding statement in the above
proposition with the same constant $c$. 
\begin{definition}
  If we wish to make precise to which Hermitian bundle $E$ we refer we shall denote the corresponding bundle as an upper-script $\mu_{0,\tau}^E$.
\end{definition}

\subsection{PU(N)-monopole equations}

In this manuscript we shall not be concerned about the analytical properties of the moduli spaces involved, such as the question of transversality. Therefore, as a matter of convenience, we shall work throughout with spaces of smooth (infinitely differentiable) sections, connections... It will be clear how to formulate a corresponding theory with Sobolev completions as is usually done in gauge theory.
\\

Let $\theta$ be a fixed smooth unitary connection in the determinant line bundle $w$. 
We shall denote by $\mathscr{A}_\theta(E)$ the space of smooth unitary connections on
$E$ which induce the fixed connection $\theta$ in $w$. This is an affine space modelled on
$\Omega^1(X;\mathfrak{su}(E))$. Here $\su(E)$ denotes the bundle of skew-adjoint
trace-free endomorphisms of $E$. Furthermore $\Gamma(X;W^+_{\data})$ denotes the space of
smooth sections of the spinor bundle $W^+_{\data}$.
We define our configuration space to be
\begin{equation*}
\mathscr{C}_{\data,\theta} := \Gamma(X;W^+_{\data}) \times \mathscr{A}_\theta(E) \ .
\end{equation*}
We denote by  $\G$ the group of unitary automorphisms of $E$ with determinant 1; it will be the `gauge group' of
our moduli problem. It acts
 in a canoncial way on sections of the spinor bundles, and as $(u,\nabla_A)
\mapsto u \nabla_A u^{-1}$ on the connections, where $u$ is a gauge
transformation and $\nabla_A$ a unitary connection. In particular it lets the induced connection in the
determinant line bundle $w$ fixed. 
The set $\mathscr{B}_{\data,\theta}$ is defined to be the configuration space up to gauge, that
is the quotient space $\conf_{\data,\theta} / \G$. 
\\

We are now able to write down the PU(N)-monopole equations associated to the
data $(\mathfrak{s},E)$ consisting of a $Spin^c$-structure $\mathfrak{s}$ and a
unitary bundle $E$ of rank $N$ on $X$. These equations read
as follows:
\begin{equation}\label{PUN-equations}
\begin{split}
\dirac_{A}^{+} \Psi & = 0 \\
 \gamma((F_{A}^+)_{0}) - \mu_{0,0}(\Psi) & = 0 \ . 
\end{split}
\end{equation}
Here $(F_{A})_{0} \in \Omega^{2}(X;\su(E))$ denotes the trace-free part of the curvature $F_{A} \in \Omega^{2}(X;\u(E))$. 

The left hand side of the above equations can be seen as a map $\mathscr{F}$ of
the configuration space $\mathscr{C}= \Gamma(X,W^+_{\data}) \times \Atheta$ to
the space
$\Gamma(X,W^-_{\data})\times \Gamma(X,\su(S^+_{\mathfrak{s}}) \tensor \su(E))$.
As such
it satisfies
the equivariance property
\[
  \mathscr{F}(u.(\Psi,{A})) = (u \times \ad_{u}) (\mathscr{F}(\Psi,{A}))
\ .
\]
Therefore it is sensible to define:
\begin{definition}
The moduli space $M_{\mathfrak{s},E,\theta}$ of $PU(N)$ monopoles is defined
to be the solution set of the equations (\ref{PUN-equations}) associated to the data $(\data)$ modulo the gauge-group $\G$:
\begin{equation*}
M_{\data,\theta} := \{ [\Psi,A] \in \bonf_{\data,\theta} | \mathscr{F}(\Psi,A) = 0 \} \ .
\end{equation*}
\end{definition}
There is an elliptic deformation complex associated to a solution $(\Psi,A)$ of the $PU(N)$ monopole equations. Its index equals minus the `expected dimension' of the moduli space. This expected dimension is computed by the Atiyah-Singer index theorem and is given by the index of the deformation operator $\delta_{A} = (-d_{A}^{*} \oplus d_{A}^{+})$ of the $PU(N)$ instantons \cite{K} plus twice the complex index of the twisted Dirac operator, 
\begin{equation}\label{index}
\begin{split}
\text{ex-dim}(M_{\data,\theta}) & = \ind(\delta_A) \oplus 2 \, \ind_\C(\dirac^+_A) \\
	& = -2 \,  \langle p_1(\su(E)),[X]\rangle - \,  (N^2-1) (b_2^+(X) - b_1(X) + 1) \\ & \ \ \ \ + \langle \text{ch}(E) e^{\frac{1}{2} c_1(S^+_\mathfrak{s})} {\hat{\text{A}}}(TX) ,
[X]
\rangle \ 
\end{split}
\end{equation}
where $p_1(\su(E))$ denotes the first Pontryagin class of the bundle $\su(E)$, $b_2^+(X)$ the maximal dimension of a subspace of $H_2(X;\R)$ on which the intersection form of $X$ is positive definite, and $\hat{\text{A}}(TX) = 1 -\frac{1}{24} \, p_1(TX)$ is the $\hat{\text{A}}$-genus of the four-manifold $X$.

\subsection{Uniform bound, compactification} \label{compactification}
The moduli space $M_{\data, \theta}$ turns out to be non-compact in general, but possesses a canonical compactification very analogue to the Uhlenbeck-compactification of instanton moduli spaces \cite{DK}. The main reason is that there is a uniform $C^0$ bound on the spinor part $\Psi$ of $PU(N)$ monopoles $[\Psi,A] \in M_{\data,\theta}$. Knowing this, the compactification is fairly standard \cite{DK}, \cite{T2}, \cite{FL3}, and therefore we will keep our exposition very brief on this point. An outline for the $PU(N)$ case can also be found in \cite{Z1}. 
\\

The $C^0$ bound is derived similarly to classical Seiberg-Witten theory \cite{KM} from the Weitzenb\"ock formula for the Dirac-operator $D_A$ by making also use of the above Proposition \ref{properness}.
\begin{prop}\label{uniform bound}
There are constants $c, K \in \R$, with $c > 0$, such that for any
monopole $[\Psi,{A}]\in M_{\data}$ we have
a $C^0$ bound:
\begin{equation}\label{apriori-bound}
	\max \abs{\Psi}^2 \, \leq \, \max \left\{ 0, K/c^2 \right\} \ . 
\end{equation}
Here the constant $K$ depends on the Riemannian metric, the fixed background
$Spin^c$ connection whereas the
constant $c$ is universal.
\end{prop}

Let $\mathfrak{s}$ be a $Spin^c$-structure on $X$ and let $E \to X$ be a
unitary bundle on $X$. We denote by $E_{-k}$ a bundle which has first Chern
class $c_1(E_{-k}) = c_1(E)$ and whose second Chern class satisfies 
\[
\langle c_2(E_{-k}), [X] \rangle = \langle c_2(E), [X] \rangle - k \ .
\]
On a four-manifold such a bundle is unique up to isomorphism.

\begin{definition}
An ideal $PU(N)$ monopole associated to the data $(\data,\theta)$ is given by a pair
$([\Psi,{A}], \bf{x})$, where $[\Psi,{A}] \in
M_{\mathfrak{s},E_{-k},\theta}$ is a $PU(N)$ monopole associated to $(\mathfrak{s},E_{-k})$ monopole, and
$\bf{x}$ is an element of the k-th symmetric power $Sym^k(X)$ of X (that is, an
unordered set of k points in X, ${\bf x} = \{ x_1, \dots, x_k\}$). The
curvature density of $([\Psi,{A}],{\bf x})$ is defined to be the measure
\[
\abs{F_A}^2 + 8 \pi^2 \sum_{x_i \in {\bf x}} \delta_{x_i} \ .
\]
The set of ideal monopoles associated to the data $(\data,\theta)$ is
\begin{equation}\label{idealmonopoles}
I M_{\data,\theta} := \coprod_{k \geq 0} M_{\mathfrak{s},E_{-k},\theta}
\times Sym^k(X) \ ,
\end{equation}
\end{definition}
which is endowed with a convenient topology \cite{DK}, \cite{T2}. Rougly speaking, in this topology a sequence in the main stratum $[\Psi_n,A_n]$ converges to a point $([\Phi,{B}],{\bf x}) \in M_{\mathfrak{s},E_{-k},\theta}
\times Sym^k(X)$ if the sequence of measures $\abs{F_{A_n}}^2 vol$ converges to the measure given by $\abs{F_{{B}}}^2 vol + 8 \pi^2 \sum_{x_i \in {\bf x}} \delta_{x_i} $, and if $\Psi_n$ converges to $\Phi$ in the complement of ${\bf x}$ in $X$. The main result is then:
\begin{theorem}(Compactness-Theorem)
The closure of $M_{\data,\theta}$ inside the space of ideal monopoles
$IM_{\data,\theta}$ is
compact.
\end{theorem}

\subsection{U(n) monopoles}

We shall denote by $\mathscr{A}(E)$ the space of smooth unitary connections on
$E$ which is an affine space modelled on
$\Omega^1(X;\mathfrak{u}(E))$. Here $\u(E)$ denotes the bundle of skew-adjoint
endomorphisms of $E$. 
We define our configuration of $U(n)$ monopoles to be the space
\begin{equation*}
\mathscr{C}_{\data} := \Gamma(X;W^+_{\data}) \times \mathscr{A}(E) \ .
\end{equation*}
We denote by  $\mathscr{G}$ the group of
{\em unitary} automorphisms of $E$ (and not just of determinant 1, as for $PU(N)$ monopoles); it is the `gauge group' of our problem. The space 
$\mathscr{B}_{\data}$ is defined to be the configuration space up to gauge, i.e. the quotient space $\conf_{\data} / \mathscr{G}$.

For a configuration $(\Psi,{A}) \in \conf_{\data}$ the $U(n)$-monopole
equations with parameter $\tau \in [0,1]$ and perturbation $\eta \in \Omega^2_+(X;i\R)$
read
\begin{equation}\label{monopole equations}
\begin{split}
\dirac^{+}_{A} \Psi & =  0 \\
 \gamma(F_{A}^+) - \mu_{0,\tau}(\Psi) & =  \gamma(\eta) \ \id   \ .
\end{split}
\end{equation}
Here $F_{{A}}$ designs the curvature of the connection ${A}$ and
$F^+_{{A}}$ its selfdual part. 

\begin{remark}
Notice that the curvature equation of (\ref{monopole equations}) splits according to the Lie algebra decomposition $\u(n) = \su(n) \oplus i \R$ into two equations:
\begin{equation*}
\begin{split}
\gamma((F_{A}^{+})_{0}) - \mu_{0,0}(\Psi) & = 0 \\
\gamma(F_{\det(A)}^{+}) - \text{\em tr} \, \mu_{0,1}(\Psi) & = n \, \gamma(\eta) \ . 
\end{split}
\end{equation*}
Therefore we see that the $U(n)$ monopole equation is in fact a coupled equation in the following sense: A solution $(A,\Psi)$ is a $PU(n)$ monopole associated to the `parameter' $\theta = \det(A)$, so that $(\det(A),\Psi)$ solves some sort of `perturbed $U(1)$ anti-self-duality equation'
\end{remark}

As above, the left hand side of the above equations can be seen as a map $\mathscr{F}_\tau$
from the configuration space $\mathscr{C}_{\data}$ to
the space
$\Gamma(X;W^-_{\data})\times \Gamma(X;\su(S^+_{\mathfrak{s}}) \tensor
\u(E)) . $
This map is equivariant with respect to the gauge group $\mathscr{G}$. The moduli space is
then defined to be the space of solutions to the monopole equations modulo
gauge:
\begin{equation*}
 M_{\data}(\tau,\eta) := \{[\Psi,{A}] \in \bonf_{\data}
| \mathscr{F}_\tau(\Psi,{A}) = (0,\gamma(\eta)) \} \ .
\end{equation*}

Again, there is an elliptic deformation complex associated to a solution 
these $U(n)$ monopole equations. The `expected dimension' of the moduli space is given by the following formula, compare formula (\ref{index}) above:
\begin{equation}\label{index}
\begin{split}
\text{ex-dim}(M_{\data}(\tau,\eta)) & =  -2 \,  \langle p_1(\su(E)),[X]\rangle - \, n^2 (b_2^+(X) - b_1(X) + 1) \\ & \ \ \ \ + \langle \text{ch}(E) e^{\frac{1}{2} c_1(S^+_\mathfrak{s})} {\hat{\text{A}}}(TX) ,
[X]
\rangle \ .
\end{split}
\end{equation}

There is also an Uhlenbeck-compactification of the moduli space $M_{\data}(\tau,\eta)$ that is stated identically to that of the $PU(N)$ monopole moduli space $M_{\data,\theta}$ in section \ref{compactification} above. In particular, the lower strata still contain bundles $E_{-k}$ with the same first Chern class as $E$. 
\\


\section{The circle-action and its fixed-point set, relations to U(n)-monopoles}
There is a circle-action on the configuration space 
modulo gauge $\bonf_{\data,\theta}$ which is induced by multiplying spinors with complex numbers of unit norm. The fixed-point set of this circle-action obviously contains elements with zero spinor component, and
the other elements are those which have a connection that splits up into the direct sum of two connections
on proper subbundles on $E$ and which have the spinor component being a section of one of these
subbundles. The latter fixed-point loci are naturally labelled by isomorphism classes of proper subbundles of
$E$. 
We shall describe a way of parametrising these fixed point loci by picking a representative vector
bundle for each isomorphism class. Next we restrict our considerations to the intersection of the
fixed-point set with the moduli space of $PU(N)$ monopoles: fixed-points with vanishing spinor are then
simply anti-selfdual $PU(N)$- connections in $E$, whereas the other fixed-point loci are fibrations of
moduli spaces of $PU(n)$-instantons in a summand $F$ of $E$ of rank $n$ over moduli spaces of
$U(N-n)$-monopoles in the complement $F^{\perp}$ of $F$ in $E$. 
\begin{remark}
In this section both configuration spaces $\conf_{\data,\theta}$ and moduli spaces $M_{\data,\theta}$ of $PU(N)$ monopoles as configuration spaces $\conf_{\s,F}$ and moduli spaces $M_{\s,F}$ of $U(n)$ monopoles appear. We remind the reader that our distinction in the notation is by the fixed connection $\theta$ in the determinant line bundle in the $PU(N)$ situation.\end{remark}

\subsection{Reductions and stabilisers of connections under the gauge group}
Here we study the stabilisers of connections $A \in \Atheta$ under
the action of the gauge group $\G$. 

An element $u$ of the gauge group $\G$ acts on the set of connections ${A}$
in $\Atheta$ by the formula 
\[
  u({A}) = {A} - (d_{{A}} u) u^{-1} \ ,
\]
where we consider $u$ as section of the vector bundle $\gl(P_E)$. The
stabiliser of the connection ${A}$ inside the gauge
group $G$ is the group of automorphisms which preserve ${A}$:
\begin{equation*}
\begin{split}
  \Gamma({A}) = & \{ u \in \G \ | \ u({A}) = {A} \} \\
   = & \{ u \in \G \ | \ d_{{A}} u = 0 \}
\end{split}
\end{equation*}
This group $\Gamma({A})$ is a finite-dimensional compact Lie-group and can  be seen as a closed Lie subgroup of $\Aut(E_x) \cong U(N)$ for any
point $x \in
X$ if $X$ ix connected. 

Note that the stabiliser always contains the centre $Z(G)$ of the
structure group $G$,
$\Gamma({A}) \supseteq Z(G)$. In our case, the centre injects as
\[
 Z(SU(N)) = \{ \lambda \ \id_E \ | \ \lambda^N = 1  \}
\]
\begin{definition}
The connection ${A} \in \Atheta$ is called {\em reducible} if the
stabiliser $\Gamma({A})$ is different from the centre $Z(SU(N))$. 
\end{definition}
Suppose now that ${A}$ is reducible and that $d_{{A}} u = 0$. We recall
a standard result for normal endomorphisms of a Hermitian vector-bundle. An
endomorphism $u$ is called normal if $u u^* = u^* u$, where $u^*$ is the adjoint
endomorphism with respect to the Hermitian structure. The following lemma is easy to prove:
\begin{lemma}\label{decomposition_of_parallel_endomorphism}
  Suppose the normal endomorphism $u \in \End(E)$ is ${A}$-parallel, $d_A u = 0$. Then
its spectrum is constant and there is a ${A}$-parallel decomposition of $E$
into subbundles
  \begin{equation*}
    E = \bigoplus_{\lambda \in \text{Spec}(u)} E_{\lambda} \ .
  \end{equation*}
Each summand $E_\lambda$ is $u$-invariant, and we
have $u|_{E_\lambda} = \lambda \ \id_{E_\lambda}$. In other words, the
$E_\lambda$ are eigen-bundles of the endomorphism $u$.
\end{lemma}
 \qed
As a corollary one obtains that a connection $A$ is reducible if and only if there is a proper
subbundle $F$ of $E$ which is $A$ - parallel. For, the latter condition clearly implies that $A$ is
reducible in our definition. On the other hand, if the stabiliser $\Gamma(A)$ is strictly bigger than
the centre $Z(SU(N))$, then by the preceding lemma there must be a $A$-parallel automorphism $u \in
\G$ which admits an eigenvalue which is not an $N^{th}$ root of one, and therefore there must be such an
$A$-parallel subbundle.

\subsection{Stabiliser of a configuration under the gauge group}

We define the stabiliser $\Gamma(\Psi,{A})$ of a configuration
$(\Psi,{A}) \in \mathscr{C}_{\data}$ to be the set 
\begin{equation*}
\begin{split}
  \Gamma((\Psi,{A})) = & \{ u \in \G \ | \ u(\Psi,{A}) =
(\Psi,{A})
\} \\
   = & \{ u \in \G \ | \ u(\Psi) = \Psi, \quad d_{{A}} u = 0 \} \
.
\end{split}
\end{equation*}

\begin{definition}
  The subset $\mathscr{C}_{\data,\theta}^* \subseteq \mathscr{C}_{\data,\theta}$
(respectively $\mathscr{C}_{\data,\theta}^{**}$) is
defined to be the set of configurations $(\Psi, {A}) \in
\mathscr{C}_{\data,\theta}$ which has zero-dimensional stabiliser (respectively
trivial stabiliser). The subset
$\mathscr{B}_{\data,\theta}^* \subseteq \mathscr{B}_{\data,\theta}$ is defined to be the
subset $\mathscr{C}_{\data,\theta}^* / \G \subseteq \mathscr{B}_{\data,\theta}$. The subset
$\mathscr{B}_{\data,\theta}^{**}$ is defined correspondingly.
\end{definition}

It is easy to see that the stabiliser of a configuration $(\Psi,{A})$, with
$\Psi$ non-vanishing and the connection ${A}$ irreducible, is trivial. 
Furthermore, the configurations $(\Psi,{A})$ with $\Psi \equiv 0$ and
irreducible connection part ${A}$ have stabilisers which are the finite
group $Z(SU(N))$. These claims, as well as the following Proposition, follow easily
from Lemma (\ref{decomposition_of_parallel_endomorphism}) in the section of
reducible connections above.

\begin{prop}\label{criterion_finite_stabiliser}
  Suppose we have given a configuration $(\Psi,{A})$ with non-vanishing 
spinor $\Psi$, and the connection of the form ${A} = {A}_1 \oplus
{A}_2$ according to a ${A}$-parallel decomposition $E=F \oplus F^\perp$,
with both ${A}_1$ and ${A}_2$ irreducible. Then its stabiliser is a
finite (and thus zero-dimensional) group.
\end{prop}

\subsection{The circle action}
We are given an $S^1$-action on the configuration space $\mathscr{C}_{\data}$
given by
the simple formula
\begin{equation*}
\begin{split}
  S^1 \times \mathscr{C}_{\data} \to & \mathscr{C}_{\data} \\
  \left(z,(\Psi,{A})\right) \mapsto & \left(z \Psi, {A}\right)
\end{split}
\end{equation*}
Now as this action commutes with the action of the gauge group $\G$, we see that
we get a well-defined action on the quotient,
\[
S^1 \times \bonf_{\data,\theta} \to \bonf_{\data,\theta} 
\]
The action is not effective. In fact if $z^N=1$ there is always a
gauge-transformation $u$ with $u(\Psi,{A}) = (z \Psi, {A})$, that is
$[\Psi,{A}] = [z \Psi, {A}]$, because $Z(SU(N)) \subseteq \Gamma(A)$. Therefore we define
\[
   r\left(z, \left[\Psi, {A}\right]\right) := \left[ z^{1/N} \Psi, {A}
\right] \ .
\]
In this formula $z^{1/N}$ is any $N^{th}$ root of $1$, the equvialence class $[ z^{1/N} \Psi, {A}]$ does
not depend on the particular choice.

\begin{remark}
Suppose we had chosen as gauge-group $\G$ the group of unitary
bundle-automorphisms of $E$ which fix the connection $\theta$ only, that is,
the larger group of unitary automorphisms with constant determinant. Then 
the same action on $\mathscr{C}_{\data,\theta}$ would
 have introduced the trivial $S^1$-action on $\mathscr{B}_{\data,\theta}$. This
justifies our choice of the gauge group $\G$ as $\Gamma(SU(E))$.
\end{remark}

\begin{lemma}
\label{fixpoints}
Suppose there exists some $z_0 \in S^1$ with $z_0\neq 1$ such that
$[z_0^{1/N}\Psi,{A}] = [\Psi, {A}]$. Then for any $z \in S^1$ we have 
  \[
    [z\Psi,{A}] = [\Psi, {A}] \ .
  \]
Thus an element $[\Psi,A]$ is a fixed-point of the circle-action if and only if it is left fixed by some
non-trivial element $z_0 \in S^1$ under the action $r$.
\end{lemma}
The proof, using Lemma \ref{decomposition_of_parallel_endomorphism}, is left to the reader.  \qed

Using the above Lemma \ref{decomposition_of_parallel_endomorphism} again we deduce the following simple criterion for fixed-points under the action $r$ above:

\begin{prop}\label{fixed points}
A configuration up to gauge $[\Psi,{A}] \in \mathscr{B}_{\data,\theta}$ is
contained in the fixed point set of the action $r$ if and only if for some (or
equivalently, for any) representative $(\Psi,{A})$ we have one of the
following (possibly both):
\begin{enumerate}
\item There is a non-trivial ${A}$-parallel orthogonal decomoposition $E=
\oplus 
E_i$ and the spinor is a section of one of $S^+_\mathfrak{s} \tensor E_i$
\item The spinor vanishes $\Psi \equiv 0$ \ .
\end{enumerate}
\end{prop}
Further down we will see that if we impose in addition the monopole equations the spinor component of a fixed point $[\Psi,A]$ will automatically lie in a proper summand $S^+ \tensor E_i$ as soon as the connection $A$ is reducible.

\subsection{The $S^1$-fixed point set inside the
configuration space modulo gauge, and parametrisations}

First we will describe the fixed-point set of the $S^1$-action inside
$\bonf_{\data,\theta}$. In the above Proposition \ref{fixed points} we saw that these are related to proper subbundles of $E$. However, two subbundles which are mapped into each other by gauge transformations, i.e. automorphisms of E, should be considered equivalent. This equivalence of subbundles might be called `ambiently isomorphic', but it is easy to see that two subbundles of $E$ are ambiently isomorphic if and only if they are isomorphic as abstract bundles. This even holds for prescribed determinant. For the further work, especially for describing the
intersection of the fixed point set with the moduli space, it turns out useful
to fix representatives $F \in [F]$ for each such isomorphism class. This 
yields to a `parametrisation' of each component of the fixed point set which is
determined by the isomorphism class $[F]$. \\

\begin{definition}\label{def_klammerf}
Let $F$ be a proper summand of the unitary bundle $E$. We define the set
$\mathscr{B}_{\data,\theta}^{[F]}$ to be the set of all elements $[\Psi,{A}]
\in \mathscr{B}_{\data,\theta}$ such that for some
representative $(\Psi,{A})$ there exists a ${A}$-parallel
decomposition $E=F \oplus F^{\perp}$ with $\Psi \in
\Gamma(X,S^+_{\mathfrak{s}} \tensor F)$.
\end{definition}
By Proposition \ref{fixed points} the space $\mathscr{B}_{\data,\theta}^{[F]}$ is contained in the
$S^1$-fixed point set $\mathscr{B}_{\data,\theta}^{S^1}$. Another subset of the fixed
point set is given by the space  $\mathscr{B}_{\data}^{\equiv 0}$ defined as the subspace of elements $[\Psi,{A}]
\in \mathscr{B}_{\data,\theta}$ with vanishing spinor component, $\Psi \equiv
0$.


The above Proposition \ref{fixed points} gives then the following description of
the fixed point set of the $S^1$ - action $r$:
\begin{prop}
The $S^1$-fixed point set $\mathscr{B}_{\data,\theta}^{S^1}$ is given as the union
\[
 \left( \bigcup_{[F] \subseteq E} \mathscr{B}_{\data,\theta}^{[F]} \right) \ \bigcup \ 
\mathscr{B}_{\data,\theta}^{\equiv 0} \ .
\]
Here the first union is taken over all isomorphism classes of
proper subbundles of $E$.
\end{prop}
We should point out as well that the different components $\bonf_{\data,\theta}^{[F]}$
may a priori intersect each other or the fixed-point component of vanishing
spinor $\bonf_{\data,\theta}^{\equiv 0}$. We would also like to remark that for rank strictly higher than $2$ we may always have
infinitely many such isomorphism classes of proper subbundles $[F]$ of $E$, even for definite intersection form. 

%

In order to have a convenient description of the set $\bonf_{\data,\theta}^{[F]}$ it seems natural to fix an actual proper subbundle $F$ for each isomorphism class $[F]$. 
Hence the following definition:

\begin{definition}
We define the configuration space relative to the splitting $E=F \oplus
F^\perp$ as the following set:
\begin{equation*}
\begin{split}
\mathscr{C}_{\data, \theta}^{F \oplus F^{\perp}} := & \left\{  (\Psi,{A}_1, 
{A}_2) \in \Gamma(S^+_\mathfrak{s}\tensor F) \times \mathscr{A}(F)
\times \mathscr{A}(F^\perp) \right| \\ &
  \left. \det({A}_1) \tensor \det({A}_2)
= \theta \right\}
\end{split}
\end{equation*}
Correspondingly, the group of unitary automorphisms with determinant 1
respecting the splitting $E= F \oplus F^\perp$ is defined to be
\begin{equation*}
\begin{split}
\mathscr{G}^0_{F \oplus F^\perp}:= \left\{ \left. (u_1,u_2) \in \Gamma(U(F))
\times \Gamma(U(F^\perp)) \right| \det(u_1) \cdot \det(u_2) = 1 
 \right\}
\end{split}
\end{equation*}
As usually, we denote the quotient by: 
\begin{equation*}
 \mathscr{B}_{\data,\theta}^{F \oplus F^\perp}:= \mathscr{C}_{\data,\theta}^{F \oplus
F^{\perp}} / \mathscr{G}^0_{F \oplus F^\perp} \ .
\end{equation*}
\end{definition}
It is then easy to see that this yields a well-defined map 
\begin{equation*}
\begin{split}
i_F: &  \mathscr{B}_{\data,\theta}^{F \oplus F^\perp} \to \mathscr{B}_{\data,\theta}^{[F]}
\\
 	& [ \Psi, {A}_1, {A}_2 ] \mapsto [\Psi, {A}_1 \oplus
{A}_2 ] \ .
\end{split}
\end{equation*}
This map is easily seen to be always surjective. However, it fails to be
injective in general. Nonetheless, on a dense subset of $\mathscr{B}_{\data}^{F
\oplus F^\perp}$ it is, as we shall show next. We will think of the map $i_F$
as a `parametrisation' of the fixed-point set component $\bonf_{\data}^{[F]}$.

As before, let us denote by
$\mathscr{B}^{* \ F \oplus F^\perp}_{\data}$ and by 
$\mathscr{B}^{* \ [F]}_{\data}$ the configurations which have
finite-dimensional
stabiliser in their groups $\G$ respectively $\mathscr{G}^0_{F \oplus F^\perp}$.
Further, we denote by $\mathscr{B}^{*,irr \ F \oplus F^\perp}_{\data,\theta}$ the
subset of $\mathscr{B}^{* \ F \oplus F^\perp}_{\data,\theta}$ consisting of
elements $[\Psi,{A}_1,{A}_2]$ with non-vanishing spinor, $\Psi \neq 0$,
and both connections ${A}_1$ and ${A}_2$ irreduible. By the
way, $(\Psi,{A}_1,{A}_2)$ has zero-dimensional stabiliser in $\G_{F
\oplus F^\perp}$ if and only if $(\Psi,{A}_1 \oplus {A}_2)$ does so in
$\G$. Now we can state the following:

\begin{prop}\label{Prop_of_iF}
Restriction of the map $i_F$ yields an injective map
\begin{equation*}
i_F: \  \mathscr{B}^{*,irr \ F \oplus F^\perp}_{\data,\theta}
\to 
\ \mathscr{B}^{* \ [F]}_{\data,\theta}
\end{equation*}
from the subset $\mathscr{B}^{*,irr \ F \oplus F^\perp}_{\data,\theta}$ of
configurations, up to gauge, with zero-dimensional stabilisers inside
$\G_{F \oplus F^\perp}$ and irreducible connections, into the fixed-point set
component $\mathscr{B}^{* \ [F]}_{\data,\theta}$ of configurations, up to gauge, with
zero-dimensional stabilisers inside $\G$.
\end{prop}
{\em Proof:}
For simplicity we note $i$ instead of $i_F$. Suppose we have elements
$[\Psi,{A}_1, {A}_2]$, $ [\Phi,{B}_1,{B}_2] \in
\mathscr{B}_{\data,\theta}^{*,irr \, F \oplus F^\perp} $ such that $i([\Psi,{A}_1,
{A}_2] = i [\Phi,{B}_1,{B}_2]$. This is equivalent to saying that
there is a gauge transformation $u \in \G$ such that 
\begin{equation*}
\begin{split}
u (\Psi) & = \Phi \\
u ({A}_1 \oplus {A}_2 ) & = {B}_1 \oplus {B}_2 \ .
\end{split}
\end{equation*}
The second equation implies that $u$ is an $({A}_1 \oplus {A}_2)
\tensor ( {B}_1 \oplus {B}_2)^*$-parallel endomorphism of $E$. Let us
write $u$ in the form 
\[
\begin{pmatrix} u_{11} & u_{12} \\ u_{21} & u_{22} \end{pmatrix}
\]
according to the splitting $E= F \oplus F^{\perp}$. Injectivity will follow if
we have $u_{12}= 0 $ and $u_{21} = 0$. It is enough to show just $u_{21}=0$,
as the other equation will follow from the fact that $u$ is unitary. We find
that the morphism $u_{11}$ is ${B}_1 \tensor {A}_1^*$-parallel,
the morphism $u_{12}$ is ${B}_1 \tensor {A}_2^*$-parallel, 
the morphism $u_{21}$ is ${B}_2 \tensor {A}_1^*$-parallel, and 
the morphism $u_{22}$ is ${B}_2 \tensor {A}_2^*$-parallel . 

Now all the connections ${A}_{i}$, ${B}_{i}$, $i=1,2$ are {unitary}
connections. Therefore the fact that, for instance, $u_{12}$ is ${B}_1
\tensor {A}_2^*$-parallel implies that the adjoint $u_{12}^*$ is ${A}_2
\tensor {B}_1^*$-parallel. As a consequence,
the endomorphism $u_{21}^* u_{21}$ of $F$ is ${A}_1 \tensor
{A}_1^*$-parallel, and 
the endomorphism $u_{12}^* u_{12}$ of $F^\perp$ is ${A}_2 \tensor
{A}_2^*$-parallel. 
By the hypothesis ${A}_1$ and ${A}_2$ are irreducible, so that the
above Lemma \ref{decomposition_of_parallel_endomorphism} implies that there are
constants $\xi, \zeta \in \C$ with 
\begin{equation*}
\begin{split}
u_{21}^* u_{21} & = \xi \ \id_{F} \\
u_{12}^* u_{12} & = \zeta \ \id_{F^\perp} \ .
\end{split}
\end{equation*}
We have to show now that under our hypothesis $\xi = 0 $ or $\zeta = 0$,
implying then that $u_{21}= 0$ respectively $u_{12}=0$. But if we had $\xi \neq
 0$, then $u_{21}$ is injective at each point $x \in X$. By the hypothesis we
get that 
$\Psi \neq 0$, and therefore we would have a non-trivial section $u_{21}
(\Psi) \in
S^+_\mathfrak{s} \tensor F^\perp$. However, we have $u(\Psi) = \Phi$, where
$\Phi$ is a section of $S^+_\mathfrak{s} \tensor F$, so that this would yield a
contradiction. Therefore $\xi = 0$ and as a consequence $u_{21}= 0$ and $u_{12}
= 0$. 
\qed 

Next we shall discuss a canonical fibering of the configuration space up to
gauge respecting the proper decomposition $E = F \oplus F^\perp$ that we have
introduced above. Let us denote now by $\G_{F}$ the group of special unitary automorphisms of the
unitary bundle $F$ on $X$, that is $\G_F = \Gamma(X,SU(F))$. So, with this
 notation, $\G_E$ is the
gauge group we have until now denoted by the letter $\G$. On the other hand, we
shall denote by $\mathscr{G}_{F}$ the group of unitary automorphisms of $F$,
that is, $\mathscr{G}_{F} = \Gamma(X,U(F))$. 
\begin{lemma}
We have an exact sequence of groups given by
\begin{equation*}
1 \to \G_{F^\perp} \stackrel{i}{\to} \G_{F \oplus F^\perp} \stackrel{j}{\to}
\mathscr{G}_{F} \to 1 \ . 
\end{equation*}
Here the morphisms are given by $i(u_2) = (\id_F, u_2)$ and $j((u_1,u_2)) :=
u_1$.
\end{lemma}
{\em Proof:} The only non-trivial point is the surjectivity of the morphism
$j$. Indeed, for a given gauge transformation $u_2 \in \mathscr{G}_{F}$ we have
to find some automorphism $u_1 \in \mathscr{G}_{F^\perp}$ such that $\det(u_1)
\cdot \det(u_2) = 1$. So we have to find an automorphism of $F^\perp$ with
prescribed determinant $\det(u_2)^{-1}$. That this is indeed possible follows
from obstruction theory \cite{S} \cite{MS}. \qed

We shall introduce some new notation now. Given a Hermitian vector bundle $F$ on $X$ we shall denote by $P_F$ its associated frame bundle, a principal bundle of structure group $U(n)$, where $n$ is the rank of $F$. Let us denote
by $\mathscr{A}_{PU}(F)$ the affine space of connections in the associated 
$PU(n)$-bundle $P_{F}\times_\pi PU(n)$, where $\pi$ is the natural projection $U(n)\to PU(n)$.
Note that in the case that $n= \text{rank}(F) = 1$ the bundle
$P_{F}\times_\pi PU(n)$ is the trivial principal bundle with structure
group the trivial group, and both $\mathscr{A}_{PU}(F)$ and
$\mathscr{A}_{PU}(F) / \G_{F}$ consist of a single point.
\begin{definition}
We shall denote by $\mathscr{B}_{F}^{PU}$ the set of all $PU(n)$-connections $A \in
\mathscr{A}_F^{PU}$ in the
unitary bundle $F$ up to the gauge group $\G_{F}$ of special unitary
automorphisms of the bundle $F$. Notice that up to specifying a connection $\vartheta \in \det(F)$ we have isomorphisms $\mathscr{A}_{F}^{PU} \equiv \mathscr{A}_{\vartheta}(F)$ and likewise $\mathscr{B}_{F}^{PU} \equiv \mathscr{B}_{\vartheta}(F)$, with the notations of Section 1. 
\end{definition}
\begin{prop}\label{map_h}
Suppose the 4-manifold $X$ is simply connected. Then we have a bijection
\begin{equation*}
\begin{split}
h: & \mathscr{B}_{\data,\theta}^{F \oplus F^\perp}  \stackrel{\cong}{\to}
\mathscr{B}_{\mathfrak{s},F} \times \mathscr{B}^{PU}_{F^\perp}
 \\
& [ \Psi,{A}_1, {A}_2 ]  \mapsto ([\Psi,{A}_1],[A_2]) \ .
\end{split}
\end{equation*}
\end{prop}
{\em Proof:}
We leave the proof of this proposition to the reader and only notice that the simply-connectedness of $X$ is used to ensure the existence of roots of $U(1)$ valued functions on $X$. In fact some (non-canonical) $(N-n)$ th root of $det(u_{1})$ will appear. For a complete proof see \cite{Z1}. \qed
\begin{remark}
  Without the assumption that $X$ is simply-connected we can still show that we
get a fibration $\mathscr{B}_{\data}^{F \oplus F^\perp} \to
\mathscr{B}_{\mathfrak{s},F}^U $ with standard fibre
$\mathscr{A}_{PU}(F^\perp)/\mathscr{G}^0_{F^\perp}$. The non-triviality of this
fibration should be encoded in $H_1(X,\Z)$. From now on, however, we shall
suppose that our 4-manifold $X$ is simply connected.
\end{remark}

\subsection{The circle-action on the moduli space of $PU(N)$-monopoles}

Until now our consideration of the $S^1$-action and its fixed point set was
inside the configuration space up to gauge, $\mathscr{B}_{\data,\theta}$. 
Obviously, the moduli space of $PU(N)$-monopoles 
$M_{\data,\theta}\subseteq \mathscr{B}_{\data,\theta}$ is invariant under the
$r$-action, $r(S^1,M_{\data}) \subseteq M_{\data}$. 
All we have found out about the circle-action on $\bonf_{\data,\theta}$ applies to the restriction of this action to the moduli space as well. However, there are more things we can say about the fixed-point set of the circle-action for this restriction. In particular, these fixed-point sets are naturally related to other moduli spaces. Obviously the intersection $\bonf_{\data,\theta}^{\equiv 0} \cap M_{\data,\theta}$ consists of anti-self-dual connections in $E$. But also the intersection $\bonf_{\data,\theta}^{[F]} \cap M_{\data,\theta}$ is parametrised by the product of the moduli space of $U(n)$-monopoles in $F$, with $n=\rk(F)$, and the moduli space of anti-self-dual $PU(N-n)$-connections in $F^\perp$, as we shall see.  

Proposition \ref{fixed points} above described the fixed-points of the circle-action on $\bonf_{\data,\theta}$, the configuration space modulo gauge. In particular the element $[\Psi,A]$ lies in $\bonf_{\data,\theta}^{[F]}$ if and only if for a representative $(\Psi,A)$ we have a (proper) $A$-parallel decomposition $E= F \oplus F^\perp$, and the spinor part $\Psi$ is a section of $S^+_\mathfrak{s} \tensor F$. This second condition becomes automatically satisfied if $(\Psi,A)$ solve the $PU(N)$ monopole equations: 

\begin{prop}\label{one_component_vanishes}
Suppose the configuration $(\Psi,{A})$ satisfies the
$PU(N)$-Seiberg-Witten-equations (\ref{PUN-equations}) associated to the data
$(\data)$. \\ Sup\-pose further that the
connection ${A}$ is reducible, and that $E = \oplus E_{i} $  is a
${A}$-parallel orthogonal decomposition into proper subbundles, and that the
base manifold $X$ is
connected. Then the spinor must be a section of one of the
bundles $S^+_\mathfrak{s} \tensor
E_i$. 
\end{prop}
{\em Proof:} Suppose the connection ${A}$ splits into two connections
${A_1}
\oplus {A_2}$ with respect to $E = E_1 \oplus E_2$. As an endomorphism of
$E$ the curvature $F_{{A}}$ splits as 
\begin{equation*}
  F_{{A}} = \begin{pmatrix} F_{{A_1}} & 0 \\ 0 & F_{{A_2}} \ 
\end{pmatrix} \ .
\end{equation*}
In other words, it is a section of $\Lambda^2(T^*X) \tensor \left(\u(E_1) \oplus
\u(E_2)\right)$. The trace-free part $F_A$ is then a section of the bundle
$\Lambda^2(T^*X) \tensor ((\u(E_1) \oplus \u(E_2))\cap \su(E))$. Therefore the
curvature-equation of the $PU(N)$-monopole-equations implies that 
\begin{equation}\label{mu-decomposition}
\mu_{0,0}(\Psi) \in \Gamma(X,\su(S^+_{\mathfrak{s}}) \tensor ((\u(E_1) \oplus \u(E_2))
\cap \su(E)))
\ .
\end{equation}
Now decompose the spinor as $\Psi= \Psi_1 + \Psi_2$, where $\Psi_i \in
\Gamma(X,S^+_{\mathfrak{s}} \tensor E_i)$. Recall that the quadratic map $\mu_{0,0}$
is defined to be
$\mu_{0,0}(\Psi) = \mu_{0,0}(\Psi,\Psi)$, where on the right we mean the bilinear map $\mu_{0,0}$.
We get
\begin{equation*}
\begin{split}
  \mu_{0,0}(\Psi,\Psi)
  	= & \mu_{0,0}(\Psi_1,\Psi_1) + \mu_{0,0}(\Psi_1, \Psi_2) + \mu_{0,0}(\Psi_2, \Psi_1) +
\mu_{0,0}(\Psi_2,\Psi_2) \ .
\end{split}
\end{equation*}

By the definition of $\mu_{0,0}$ and by the above
equation (\ref{mu-decomposition}) we see that $\mu_{0,0}(\Psi_1,\Psi_2) =
\mu_{0,0}(\Psi_2,\Psi_1) = 0$. Now from the fact that the bilinear map $\mu_{0,0}$ is
`without zero-divisors' by the above Proposition \ref{properness} we see that
in each fibre $\Psi_1= 0 $ or $\Psi_2=0$. Suppose we have $\Psi_1(x_0) \neq 0$
for some point $x_0 \in X$. As $\Psi$ is continuous we must have
$\Psi \neq 0$ for all $x$ in some neighbourhood $U$ of $x_0$. Therefore
$\Psi_2 \equiv 0$ on $U$. However, the Dirac equation $\dirac_{{A}}^{+} \Psi = 0$
implies that $\dirac_{{A}_1}^{+} \Psi_1 = 0$ and that $\dirac_{{A}_2}^{+} \Psi_2 = 0$,
where the Dirac operator $D_{{A}_i}^{+}: \Gamma(X,S^+_{\mathfrak{s}} \tensor
E_i) \to
\Gamma(X,S^-_{\mathfrak{s}} \tensor E_i)$ is defined to be the composition of
$\nabla_{{A_i},B}: \Gamma(X, S^+_{\mathfrak{s}} \tensor E_i) \to \Omega^1(X,
S^+_{\mathfrak{s}} \tensor
E_i)$ with the Clifford-map $\gamma: T^*X \tensor (S^+_{\mathfrak{s}} \tensor
E_i) \to (S^-_{\mathfrak{s}}
\tensor E_i)$. But for each of these Dirac operators there is a unique
continuation theorem for elements in its kernel by Aronaszajin's Theorem
\cite{A}. Therefore, as $\Psi_2 \equiv 0 $ on $U$, it must vanish identically on
$X$. The general case follows easily by iterating the same argument. \qed

\begin{remark}
The $S^1$-action extends naturally to the Uhlenbeck-com\-pacti\-fication
\[
\bar{M}_{\data,\theta} \ \subset \ I M_{\data,\theta} =  \coprod_{k \geq 0}
M_{\data_{-k},\theta} \times Sym^k(X) \ 
\]
(with its above-defined topology). 
\end{remark}
Another important result is the following finiteness property of the fixed-point locus inside the (compactified) moduli space:
%
%

\begin{prop}
Given the data $(\data)$ the respective moduli spaces $M_{\s,E_{-k},\theta}$, $k \geq 0$, occuring in the definition of the Uhlenbeck compactification of $M_{\data,\theta}$,  intersect the respective fixed point loci $\mathscr{B}_{\s,E_{-k},\theta}^{[F]}$ only in finitely many isomorphism classes of proper summands $[F]$ of $E$.
\end{prop}

{\em Proof:}
We will show that if $[\Psi,{A}] \in M_{\data} \cap
\mathscr{B}_{\data}^{[F]}$,
then $c_1^{\R}(F)$ lies in a bounded set within $H^2(X,\R)$, and 
$\langle c_2(F),[X] \rangle \in \Z$ is bounded also. As $c_1^{\R}(F)$ is in the image of the morphism $H^2(X,\Z) \to
H^2(X,\R)$, it will follow that $c_1(F)$ lies in a finite set.
The conclusion is
then that only finitely many pairs $(c_1,c_2) \in H^2(X,\Z) \times H^4(X,\Z)$
can occur as first and second Chern-class of $F$. But on a closed oriented
4-manifold unitary bundles are classified, up to isomorphism, by their first
and second Chern class.

Recall the Chern-Weil formulae for the image of the first and second Chern class
inside $H^*(X,\R) \cong H^*_{dR}(X)$:
\begin{equation}\label{Chern-Weil}
\begin{split}
  c_1^{\R}(E) = & \ \frac{-1}{2\pi i} \left[ \tr F_{{A}} \right] \ , \\
  c_2^{\R}(E) = & \ \frac{-1}{4\pi^2} \left[ \frac{1}{2} \left( \tr{F_{{A}}}
\wedge
\tr{F_{{A}}} - \tr ( F_{{A}} \wedge F_{{A}}) \right) \right] \\
    = & \frac{1}{2} \langle c_1(E)^2, [X] \rangle 
 	+ \frac{1}{8\pi^2} \left( \norm{F_{{A}}^-}_{L^2(X)}^2 
		- \norm{F_{{A}}^+}_{L^2(X)}^2 \right) \ .
\end{split}
\end{equation}

The vector space $H^{2}(X;\R)$ is isomorphic to the space of harmonic 2-forms $\mathscr{H}^{2}(X,g)$. For each class $[\omega] \in H^{2}(X;\R)$ its harmonic representative $\omega_{g}$ is minimising the $L^{2}$ norm among all representatives of the same class. So if we give the space $H^{2}(X;\R)$ the inner product via its identification with $\mathscr{H}^{2}(X,g)$ we see that a subset of $H^{2}(X;\R)$ is bounded if the $L^{2}$ norms of a set of representative forms is bounded. Thus to bound the classes $c_1(F) = [\tr
F_{{A}_1}]$ it is enough to bound the norms 
\[
  \norm{\tr F_{{A}_1}}_{L^2(X)} \ 
\]
with $A_{1}$ a connection on the summand $F$ of $E$.

Now by the assumption that $[\Psi,{A}] \in M_{\data} \cap
\mathscr{B}_{\data}^{[F]}$ we have a connection ${A}$ on $E$ that reduces to ${A}_1 \oplus
{A}_2$ according to the splitting $E = F \oplus F^{\perp}$. We therefore get the decomposition
\[
  F_{{A}} = \begin{pmatrix} F_{{A}_1} & 0 \\ 0 & F_{{A}_2}
\end{pmatrix} \ .
\]
In particular, we have bounds
\begin{equation*}
\begin{split}
\norm{\tr F_{{A}_1}}_{L^2(X)} & \leq \norm{\tr F_{{A}}}_{L^2(X)} \ , \\
\norm{F_{A_{1}}^{+}}_{L^{2}(X)} & \leq \norm{F_{A}^{+}}_{L^{2}(X)} \, \\
\norm{F_{A_{1}}^{-}}_{L^{2}(X)} & \leq \norm{F_{A}^{-}}_{L^{2}(X)} \ .
\end{split}
\end{equation*}

However, by the $PU(N)$ monopole equations (\ref{PUN-equations}), the a-prioi bound (\ref{apriori-bound}), and the Chern-Weil formula for $c_{2}$ we see that the quantities on the right hand sides are bounded given by expressions that depend on the metric $g$, the
$Spin^c$ connection $B$ and the connection $\theta$ in the determinant line
bundle, as well as on some constants related to the $\mu$-map and $\gamma$, so these are uniformly bounded on $M_{\data,\theta}$. 

By the Chern-Weil formulae we therefore see that $c_1(F)$
is uniformly bounded in $H^2(X,\R)$, and therefore $ \abs{\langle c_2(F), [X] \rangle }$ is likewise bounded. As the apriori-bound on the spinor (\ref{apriori-bound}) does not depend on the second Chern class of $E$ it follows that there is a corresponding statement for the lower strata of the Uhlenbeck compactification.

\qed

%
In the sequal we shall denote by $M_{\data}^{S^1}$ the intersection
$\mathscr{B}_{\data}^{S^1} \cap M_{\data}$, as well as by
$M_{\data}^{[F]}$ the intersection of $\mathscr{B}_{\data}^{[F]} \cap
M_{\data}$. Also, $M_{\data}^*$ shall denote the
intersection of $\mathscr{B}_{\data}^*$ with $M_{\data}$, and 
$M^*{}_{\data}^{S^1}$, $M_{\data}^*{}^{[F]}$ the
respective intersections with the fixed point set and the given fixed point set
component.

\subsection{Monopole equations for configurations mapping to the fixed point
set}
Above we have pointed out that for describing the component of the fixed point
set $\mathscr{B}_{\data}^{[F]}$ determined by the isomorphism class of a proper
subbundle $[F]$ of $E$, it is useful to keep a representative $F$ fixed. We did
then describe the component $\mathscr{B}_{\data}^{[F]}$ as the image via
$i_F$ of the space $\mathscr{B}_{\data}^{F \oplus F^\perp}$ which is easier
to handle with.
It will turn out that this way we also get a convenient description of
$M_{\data}^{[F]}$, which we define to be the intersection of $\bonf_{\data}^{[F]}$ with the moduli space $M_{\data}$. 

Let us write down explicitly the
monopole equations which are satisfied by a representative $(\Psi,{A})$
having the property that there is a ${A}$-parallel decomposition of $E$
into $F \oplus F^\perp$, with $\Psi$ a section of $S_{\mathfrak{s}}^+ \tensor
F$, and ${A}$ splitting as ${A}_1 \oplus {A}_2$. Recall that $\det({A}) = \det(A_{1}) \tensor \det(A_2)$ is the fixed connection $\theta$ in
the determinant line bundle $\det(E)$. We then have 
\[
  (F_{A})_{0} = F_{A} - \, \frac{1}{N} \, \tr( F_{A}) \, \id_E = \begin{pmatrix} F_{{A}_1} -\frac{1}{N} \, F_\theta & 0 \\ 0 & F_{{A}_2} - \frac{1}{N} \, F_\theta
\end{pmatrix} \   ,
\]
according to the splitting $E = F \oplus F^\perp$, and also

\begin{equation*}
\begin{split}
\mu_{0,0}(\Psi) & = \begin{pmatrix} \mu^{F}_{0,1}(\Psi) - \ \frac{1}{N} \ \tr \ 
\mu_{0,1}^{F}(\Psi) \ \id_{F} & 0 \\ 
0 & -\frac{1}{N} \ \tr \ \mu_{0,1}^{F}(\Psi) \ \id_{F^\perp} \end{pmatrix} \\
	& = \begin{pmatrix} \mu^{F}_{0,1-\frac{n}{N}}(\Psi) & 0 \\ 
0 & -\frac{1}{N} \ \tr \ \mu_{0,1}^{F}(\Psi) \ \id_{F^\perp} \end{pmatrix} \ . \end{split}  \end{equation*}
The $PU(N)$-monopole equations (\ref{PUN-equations}) for the pair $(\Psi, {A}_1 \oplus A_2)$ then read
\begin{equation}\label{split_monopoles prov}
\begin{split}
\dirac_{A_{1}}^{+} \Psi & =  0 \\
 \gamma(F_{{A}_1}^+) -
\mu_{0,1-\frac{n}{N}}^F(\Psi) & = \frac{1}{N} \, \gamma ( F_\theta^+) \ \id_{F} \\
 \gamma(F_{{A}_2}^+ ) +
\frac{1}{N} \ \tr \ \mu_{0,1}^F(\Psi) \  \id_{F^\perp} & =  \frac{1}{N}\, \gamma
( F_\theta^+) \ \id_{F^\perp} \ . \\
\end{split}
\end{equation}
Here the terms in the second equation are sections of the bundle $\su(S_\mathfrak{s}^+)\tensor_\R \u(F)$ and the terms in the third equation are sections of $\su(S_\mathfrak{s}^+) \tensor_\R \u(F^\perp)$. 
There is Lie algebra decompositions $\u(F) = \su(F) \oplus i \R$ and correspondingly for $\u(F^\perp)$. It turns out that the `$i \R$' component of the second and the third equation are equivalent. Indeed, taking the trace (with respect to the factor $\u(F)$ in 
$\su(S_\mathfrak{s}^+)\tensor_\R \u(F)$, and correspondingly for $\u(F^\perp)$) of the second and the third equation, and using the fact that 
\[
\tr(F_{A_1}) + \tr(F_{A_2}) = F_\theta \ ,
\]
this follows from a simple computation. Therefore the system of equations (\ref{split_monopoles prov}) above is equivalent to the same system where we take as the third equation only the component of $\su(F^\perp)$ according to $\u(F^\perp) = \su(F^\perp) \oplus i \R $. Thus the $PU(N)$ monopole equations are therefore equivalent to 
\begin{equation}\label{split monopoles}
\begin{split}
\dirac_{A_{1}}^{+} \Psi & =  0 \\
 \gamma (F_{{A}_1}^+ ) -
\mu_{0,1-\frac{n}{N}}^F(\Psi) & = \frac{1}{N} \, \gamma ( F_\theta^+ \ \id_{F}) \\
(F_{A_2}^+)_{0}   & =  0 \ . \\
\end{split}
\end{equation}
The first two equations of (\ref{split monopoles}) are $U(n)$- monopole equations for $(\Psi,A_1)$ with parameters $\tau=1-\frac{n}{N} \in [0,1]$ and self-dual 2-form $\eta = \frac{1}{N} \,  F_{\theta}^{+}$, and the third equation is the anit-self-duality equation for the $PU(N-n)$ connection $A_2$. We shall denote by $M_{F}^{asd}$ the moduli space of anti-self-dual $PU(n)$ - connections in $F$ which is defined to be the space of $PU(n)$- connections $A \in \mathscr{A}^{PU}(F)$ in $F$ which satisfy the equations $F_A^+ = 0$, quotiented by the action of the gauge-group $\G_{F}$ of special unitary automorphisms of $F$. Equivalently, if we think of $A$ as a unitary connection under an isomorphism $\mathscr{A}^{PU}(F) \cong \mathscr{A}_{\vartheta}(F)$ specified by a fixed connection $\vartheta$ in the determinant line bundle $\det(F)$, the (projective) anti-selfduality equation for $A$ becomes $(F_{A}^{+})_{0} = 0$. 

We summarise this computation in the following:
\begin{prop}\label{decoupled equations}
Suppose the configuration $(\Psi,A) \in \conf_{\data}$ has reducible connection part $A = A_1 \oplus A_2$ according to $E = F \oplus F^\perp$, and that the spinor part $\Psi$ is a section of $S^+ \tensor F$ (compare proposition \ref{one_component_vanishes}). Then the $PU(N)$ monopole equations for $(\Psi,A)$ are equivalent to the system (\ref{split monopoles}). In particular, the configuration $(\Psi,A_{1})$ represents a $U(n)$ monopole in the moduli space $M_{\s,F}(1-\frac{n}{N},\frac{1}{N} F_{\theta}^{+})$, and the connection $A_{2}$ represents an instanton in the moduli space $M^{asd}_{F^{\perp}}$.
\end{prop}


\begin{definition}
We shall denote by $M_{\data}^{F \oplus F^\perp} \subseteq
\mathscr{B}_{\data}^{F \oplus F^\perp}$ the moduli space space of solutions
$(\Psi,{A}_1,{A}_2) $ to the above equations
(\ref{split monopoles}) modulo the gauge group $\G_{F \oplus F^\perp}$.
As usually, we denote by $M_{\data}^{*\ F \oplus F^\perp}$ the subspace
of those elements whose representatives have zero-dimensional stabiliser.
\end{definition}

\begin{prop}\label{iF_modulispace}
The map $i_F:\mathscr{B}_{\data}^{F\oplus F^\perp} \to
\mathscr{B}_{\data}^{[F]}$ maps the moduli space $M_{\data}^{F \oplus
F^\perp}$ onto the fixed point set component $M_{\data}^{[F]}$ inside
the moduli space. It maps the set $M_{\data}^*{}^{F \oplus F^\perp}$
bijectively onto $M^*{}_{\data}^{[F]}$.

\end{prop}
{\em Proof:} The fact that the map is onto is an immediate consequence of the above Proposition \ref{decoupled equations} and the definition of $M_{\data}^{[F]}$. 
For the remaining claim we will show that we can apply the above Proposition \ref{Prop_of_iF}. First, we shall observe that
 if $[\Psi,{A}_1,{A}_2]$ belongs to $M_{\data}^{*\ F \oplus
F^\perp}$, then the connections ${A}_1$ and ${A}_2$ are indeed
irreducible. Obviously $A_2$ has to be irreducible, but suppose ${A}_1$ were reducible. We would have a ${A}_1$ - parallel orthogonal decomposition $F = F_1 \oplus
F_2$, with ${A}_1$ splitting accordingly, ${A}_1 = {A}_{11} \oplus
{A}_{12}$. 

Let us write $\Psi=\Psi_1 + \Psi_2 \in
\Gamma(X,S^+_{\mathfrak{s}} \tensor (F_1 \oplus F_2))$ for the corresponding decomposition of
the spinor. We claim that either $\Psi_1 = 0$ or $\Psi_2 = 0$.
In fact, $(\Psi,A_1)$ solves the first two of the equations (\ref{split monopoles}). The map $\mu_{0,\tau}^F$ is `without zero-divisors' by the above Proposition \ref{properness}. With this fact the conclusion follows exactly like in the proof of Proposition \ref{one_component_vanishes}. But then the configuration
$(\Psi,{A}_{11}\oplus {A}_{12},{A}_2)$ must have
positive-dimensional stabiliser inside $\G_{F \oplus F^\perp}$, and the element 
$[\Psi,{A}_1,{A}_2]$ would not belong to $M_{\data}^{*\ F
\oplus
F^\perp}$. Therefore $[\Psi,{A}_1,{A}_2]$ belongs to the set
$\mathscr{B}_{\data}^{*,irr \ F \oplus
F^\perp}$ and we can apply Proposition $\ref{Prop_of_iF}$ for getting
injectivity. Furthermore it is easy to see that the parametrisation $i_F$
maps $M_{\data}^*{}^{F \oplus F^\perp}$
onto $M^*{}_{\data}^{[F]}$.
\qed

\begin{prop}\label{moduli_space_product}
  Restricting the bijection $h$ of Proposition \ref{map_h} above to the moduli
space $M_{\data,\theta}^{F \oplus F^\perp}$ we get an induced bijection
\begin{equation*}
h|_{M} : M_{\data,\theta}^{F \oplus F^\perp} \stackrel{\cong}{\to}
\
M_{\mathfrak{s},F}(1-n/N,1/N\, F_{\theta}^{+}) \times M^{asd}_{F^\perp} 
\end{equation*}
Together with the map $i_F$ we thus get the parametrisation of the fixed point
set component $M_{\data}^{[F]}$ of the moduli space as the
product of
a moduli space of $U(n)$-monopoles with the moduli space of
$ASD-PU(N-n)$-connections. In particular, for the irreducible parts we get a
bijection
\begin{equation*}
i_F \circ h|_M^{-1} : M_{\mathfrak{s},F}^{*}(1-n/N,1/N\, F_{\theta}^{+}) \times
M^{* \ asd}_{\ F^\perp}  \stackrel{\cong}{\to}
M_{\data,\theta}^{* \ [F]} \ .
\end{equation*}
\end{prop}
This follows from Proposition \ref{iF_modulispace} and Proposition \ref{map_h},
where it is easily checked that $h|_M^{-1}$ maps the `irreducibles'
$M_{\dataf}^{*} \times M_{F^\perp}^{* \, asd}$ onto the
corresponding `irreducibles' $M_{\data,\theta}^{* \ F\oplus F^\perp}$. 
\qed 

The whole discussion is now summarised in 

\begin{theorem}\label{S1fixedpointset_main_text}
The fixed point set under the above circle-action $r$ on the moduli space $M_{\data,\theta}$ of $PU(N)$  monopoles is given as the union of the moduli space $M_E^{asd}$ of
anti-self-dual $PU(N)$ connections in $E$ and a finite union  
\[
  \bigcup_{[F] \subseteq E} M_{\data}^{[F]} 
\]
 of components $M_{\data}^{[F]}$ indexed by a finite number of isomorphism classes $[F]$  of
proper subbundles of $E$. The spaces $M_{\data}^{[F]}$ are given as
follows: An element $[\Psi,{A}]$ belongs to $M_{\data}^{[F]}$ if 
for each representative $F \in [F]$ there is a
representative $(\Psi,{A}) \in [\Psi,{A}]$ such that $F$ is an
${A}$-invariant proper subbundle of $E$ and the spinor $\Psi$ is a section of the
proper subbundle $S^+_\mathfrak{s} \tensor F$ of $W^+_{\data}=S^+_\mathfrak{s} \tensor E$.

 Furthermore, if $X$ is simply connected, then there is a parametrisation of
this
space $M_{\data,\theta}^{[F]}$ as the product
\[
  M_{\mathfrak{s},F}(1-n/N,1/N\, F_{\theta}^{+}) \times M_{F^\perp}^{asd} \to
M_{\data,\theta}^{[F]} \ .
\]
This map is a surjection and is a bijection between the open and dense subsets of elements with zero-dimensional stabiliser in the corresponding moduli spaces, 
\[
  M_{\mathfrak{s},F}^{*}(1-n/N,1/N\, F_{\theta}^{+}) \times M_{F^\perp}^{*\ asd} \stackrel{\cong}{\to}
M_{\data,\theta}^{*\ [F]} \ .
\]
\end{theorem}
\qed
We observe that there is a corresponding statement if we take the whole Uhlenbeck-compactification of $M_{\data,\theta}$ into account. 
\section{What to expect from the cobordism program} 

The heuristical idea of the cobordism program is that the $S^{1}$-quotient of the complement of the fixed point locus $M_{\data,\theta} \setminus M_{\data,\theta}^{S^{1}}$ yields an oriented cobordism between 
\begin{enumerate}
\item a projective bundle over the instanton moduli space $M^{asd}_{E}$, the fibre over an instanton $A$ being the projectivisation of the kernel of the dirac operator $\dirac_{A}^{+}$, and
\item projective bundles over the moduli spaces $M_{\s,F} \times M^{asd}_{F^{\perp}}$, corresponding to a parallel decompositions $E = F \oplus F^{\perp}$, the fibres being projectivisation of the (complex) normal bundles that can be described by local models around these $S^{1}$ fixed point spaces. 
\end{enumerate}
Then there are extension of the Donaldson-$\mu$-classes (see \cite[Section 5 and 9]{DK} for instance), lifted to the projective bundle occuring in (1) (and multiplied with a power of the first Chern class of the bundle of kernels of the Dirac operator), to the whole $S^{1}$-quotient. Therefore, the higher rank instanton invariants, obtained by evaluating the $\mu$-classes on the `fundamental cylce' given by the moduli space $M^{asd}_{E}$, should be expressible by corresponding evaluations involving the spaces $M_{\s,F} \times M^{asd}_{\s,F^{\perp}}$. 

The components of the second type involving $M_{\s,F}$ for $rk(F) = 1$ are contributions to this evaluation that are expected to involve the Seiberg-Witten invariants of the $Spin^{c}$ structure $\s \tensor F$. The corresponding generalisation of Witten's conjecture would follow then by induction on the rank $N$ if the the moduli spaces $M_{\s,F}$ for $rk(F) > 1$ contributed trivially to a cobordism-formula for the higher rank invariant to which we made allusion. We present two arguments why such contributions can be expected to vanish. The second, by a consideration on K\"ahler surfaces, is obtained in the next section. The first, admittedly less convincing, shall be outlined in the remainder of this section.
\\

First, let us recall that the $U(n)$ monopole moduli spaces $M_{\s,F}(\tau)$ admit an Uhlenbeck compactification for any $\tau \in [0,1]$ for we still obtain an apriori-bound if $n > 1$. Let us assume we have a space of perturbations $\mathscr{P}$ such that the irredubible part of the `parametrised moduli space', the zero-locus of the resulting map
\begin{equation*}
	\bonf_{\s,F}^{**} \times \mathscr{P} \times [0,1] \to \Gamma(X;S_{\s}^{-} \tensor F \oplus \Lambda^{2}_{+}(\u(F)) \ ,
\end{equation*}
with $[0,1]$ being the parameter space for $\tau$, is cut out transversally. The existence of such a parameter space, compatible with some kind of Uhlenbeck compactification for the moduli space, is not speculative, see for instance \cite{FL4} for the case $n = 2$ or \cite{Z} for general $n$ (but less general with respect to Uhlenbeck-compactification). We may furthermore assume that the restriction to any $\tau \in [0,1]$ is transversal. 

Let us denote the perturbation parameter by $p \in \mathscr{P}$, and let us further assume that the restriction of the parametrised moduli space to $(\tau, p)$, denoted by $M_{\s,F}(\tau,p)$, is already compact in the main stratum for all $\tau \in [0,1]$ and all `small enough p'. We shall furthermore assume that $M_{\s,F}(\tau,p) = M^{**}_{\s,F}(\tau,p)$, that is, that no reducibles occur for generic parameters $p$ - this being an admittedly speculative assumption that we don't expect to hold in general, even if $b_{2}^{+}(X)>0$. Then for generic parameter $p$ the moduli space $M_{\s,F}(\tau,p)$ is a smooth closed manifold of the expected dimension for generic perturbation $p$. \\

By the usual argument the closed manifolds  $M_{\s,F}(\tau,p)$ and $M_{\s,F}(0,p)$ are then cobordant. However, it is a rather easy observation that the moduli spaces $M_{\s,F}(0,p)$ are generically empty. In fact, for parameter $\tau = 0$ the map $\mu_{0,\tau}$ is traceless, so the trace of the curvature equation of (\ref{monopole equations}) for $(A,\Psi)$ becomes:
\begin{equation}\label{abelian asd}
	F^{+}_{\det(A)} =  \text{\em pr}_{\Lambda^{2}_{+}\tensor i \R} (p) . 
\end{equation}
This is a perturbed abelian $ASD$ equation. As the derivative of the `map' $F^{+}$ on $U(1)$ connections is given by $ \, d^{+} : \Omega^{1}(X;i \R) \to \Omega^{2}_{+}(X;i \R)$, and as the cokernel of this map has dimension $b_{2}^{+}(X)$ we see that for $b_{2}^{+}(X) > 0$ the space of solutions $(A,\Psi)$ to (\ref{abelian asd}) is empty for generic parameter $p \in \mathscr{P}$. 

\section{$U(n)$ moduli spaces on K\"ahler surfaces}
In classical Seiberg-Witten theory K\"ahler surfaces are of a significant importance. Indeed, they provided the first examples of 4-manifolds with non-trivial Seiberg-Witten invariants \cite{W}. This was generalised to symplectic manifolds \cite{Ta}. All other non-vanishing results known to the author are derived from these manifolds by various kinds of glueing results for the Seiberg-Witten invariants \cite{Ta}, \cite{Fr}.

As the $U(n)$ monopole equations are a generalisation of the classical Seiberg-Witten equations it is therefore most natural to study the $U(n)$ monopole moduli spaces for K\"ahler surfaces. Whereas the analysis of the $U(n)$ monopole equations on K\"ahler surfaces is very analogous to the classical situation the final conclusion is in sharp contrast to the classical situation. Indeed, we will show in Corollary \ref{final conclusion} that if we perturb the monopole moduli space on a K\"ahler surface with a non-vanishing holomorphic 2-form then the associated moduli space is empty. 

Non-abelian monopoles on K\"ahler surfaces have also been studied by Teleman \cite{T}, Okonek and Teleman \cite{OT3} and by Bradlow and Garcia-Prada \cite{BG}, but with a rather complex geometric motivation. Corollary \ref{final conclusion} seems to appear here for the first time.

\subsection{The $U(n)$ - monopole equations on K\"ahler surfaces}
We will quickly recall now the canonical $Spin^c-$ structure on an almost complex surface. The additional
condition of $X$ being K\"ahler implies that there is a canonical $Spin^c$ connection induced by the
Levi-Civita connection. This will be our fixed back-ground $Spin^c$ connection and it is then simple to
determine the Dirac-operator associated to this fixed connection and a $U(n)-$ connection in a Hermitian
bundle $E$. We will then write down the $U(n)$ monopole equations in this particular setting.
\\

Suppose we have an almost complex structure $J: TX \to TX$ on the closed,
oriented Riemannian 4-manifold $X$ which is isometric. The associated K\"ahler
form
$\omega$ is defined by the formula
\[
   \omega_g(v,w):= g(Jv,w) \ .
\]
This is an anti-symmetric form of type $(1,1)$ when extended to the
complexification $TX^\C := TX
\tensor_\R \C$. It is a fundamental fact that the complexification of the bundle
of self-dual two forms is given by 
\[
  \Lambda^2_+ \tensor \C = \C \omega_g \oplus \Lambda^{2,0} \oplus \Lambda^{0,2}
\ .
\]
Let $e(u)$ denotes exterior multiplication with the form $u \in
\Lambda(T^*X^\C)$
and $e^*(u)$ its adjoint with respect to the inner product induced by the
Riemannian metric.

There is a canonical $Spin^c$-structure 
associated to an
almost-complex structure $J$ on $X$ \cite{Hi}. We shall denote it by
$\mathfrak{c}$. The
spinor bundles are defined to be 
\begin{equation*}
\begin{split}
  S^+_{\mathfrak{c}} := & \Lambda^{0,0}(X) \oplus \Lambda^{0,2}(X) \ ,\\
  S^-_{\mathfrak{c}} := & \Lambda^{0,1}(X) \ ,
\end{split}
\end{equation*}
and the Clifford multiplication is given by 
\begin{equation*}
\begin{split}
 \gamma  & :  \Lambda^1(T^*X) \to
\Hom_{\C}(S^+_{\mathfrak{c}},S^-_{\mathfrak{c}}) \\
   & u \mapsto \sqrt{2} (e(u^{0,1}) - e^*(u^{0,1}))  \ .
\end{split}
\end{equation*}
The induced isomorphism
\[
\gamma : \Lambda^2_+(X) \tensor \C \to \mathfrak{sl}(S^+_{\mathfrak{c}})
\]
is then seen to be given by the formula
\begin{equation}\label{cliffordmultiplication}
\gamma(\eta^{1,1} + \eta^{2,0} + \eta^{0,2}) = 
4 \begin{pmatrix} -i \Lambda_g(\eta^{1,1}) & - * (\eta^{2,0} \wedge \_ \ )
\\ \eta^{0,2} & i \Lambda_g(\eta^{1,1}) \end{pmatrix} \ .
\end{equation}
Here we use the commonly used convention to denote contraction with $\omega_g$,
that is $e^*(\omega_g)$, by the symbol $\Lambda_g$.

Now suppose that $X$ is a K\"ahler surface. This means that first the almost
complex structure $J$ is integrable to a complex structure, and second that the
K\"ahler form $\omega_g$ is closed, $d \omega_g = 0$. The condition of
closedness implies (cf. \cite{KN2}, p. 148) that the the almost complex
structure $J$ is parallel with respect to the Levi-Civita-connection
$\nabla_g$. As a consequence, the splittings 
\[
\Lambda^k(X)\tensor \C =
\oplus_{p+q=k}\Lambda^{p,q}(X)
\]
are $\nabla_g$-parallel, where we also denote by $\nabla_g$ the connection induced by the Levi-Civita
connection on all exteriour powers of $T^*X$. The canonical $Spin^c$-connection is now simply given by the
the connection $\nabla_g$ in the bundles $\Lambda^{0,0},
\Lambda^{0,1} $ and $\Lambda^{0,2}$.

%

Let $E$ be a Hermitian vector bundle on $X$, and further $\nabla_{{A}}$
a unitary connection on $E$. We shall use the notation convention
$\Lambda^{p,q}(E) := \Lambda^{p,q}(X) \tensor E$, and by $\Omega^{p,q}(E)$ we
shall denote the space of sections of the latter bundle,  $\Omega^{p,q}(E) =
\Gamma(\Lambda^{p,q}(E))$ .
\begin{definition}
  The operator $\bar{\partial}_{{A}} : \Omega^{p,q}( E)
\to \Omega^{p,q+1}( E)$ is defined to be the composition of 
$d_{{A}} : \Omega^{p+q}(E) \to
\Omega^{p+q+1}(E)$, the extension of the exteriour derivative to forms with values in $E$ by means of the connection $\nabla_{A}$,
 with the bundle projection
$\Lambda^{p+q+1}(E) \to
\Lambda^{p,q+1}( E)$.
\end{definition}

The Dirac operator associated to the canonical
$Spin^c$-connection
$\nabla_g$ in the canonical $Spin^c$-structure $\mathfrak{s}_c$ and the unitary
connection ${A}$ in the Hermitian bundle $E$ 
is expressible in terms of the above operator $\bar{\partial}_{{A}}$ and
its formal $L^2$-adjoint $\bar{\partial}_{{A}}^*$ as follows:
\begin{equation} \label{dirackaehler}
  \dirac_{{A}}^{+} = \sqrt{2} \left( \bar{\partial}_{{A}} +
\bar{\partial}_{{A}}^* \right) \ .
\end{equation}
This is a well-known fact in the case $n=1$ \cite{Hi}. The proof of the general case follows along the same
lines. In particular, the proof given in the lecture notes \cite{T3} is directly applicable to our
situation. \qed 

%

We will now study the $U(n)$ monopoles associated to the data
$(\mathfrak{c},E)$ with spinor bundles $W_{\mathfrak{c},E}^\pm = S^\pm \tensor
E$. Note that, up to tensoring $E$ with a line bundle, we can always
assume that general data $(\mathfrak{s},E)$ is of the particular form
$(\mathfrak{c},E)$. Now according to the isomorphism
$W^+_{\mathfrak{c},E} \cong \Lambda^{0,0}(E) \oplus \Lambda^{0,2}(E)$ a
spinor $\Psi \in \Gamma(X;W^+_{\mathfrak{c},E})$ can be written as
$\Psi=(\alpha,\beta)$ with $\alpha \in \Omega^{0,0}(X;E)$ a section of $E$ and
$\beta \in \Omega^{0,2}(X;E)$ a 2-form of type $(0,2)$ with values in $E$.
We introduce the following notations. We denote by $ ^{-}:
\Lambda^{p,q}(E) \to \Lambda^{q,p}(E^*)$ the conjugate linear isomorphism which
is the tensor product of complex conjugation on the forms and the conjugate
linear isomorphism specified by the hermitian structure on the bundle $E$. We
denote by $^* : \Lambda^{p,q}(E) \to \Hom(\Lambda^{p,q}(E),\C)$ the conjugate
linear isomorphism specified by the Hermitian structure on $\Lambda^{p,q}(E)$.
For an endomorphism $f \in \End(E)$ we denote $\{ f \}_\tau := (f)_0 +
\frac{\tau}{n} \tr(f) \id_E $, where $(f)_0$ denotes the trace-free part of
$f$. Thus we simply have $\{f\}_1 = f$. With this said we can write
$\mu_{0,\tau}(\Psi)$ according to the above isomorphism as 

\begin{equation}\label{mukaehler}
  \mu_{0,\tau}(\Psi) = \begin{pmatrix}
  \frac{1}{2} \left( \{\alpha \alpha^*\}_\tau - \{ * \beta \wedge \bar{\beta}
\}_\tau \right)  & \{\alpha \beta^*\}_\tau \\
\{\beta \alpha^*\}_\tau & \frac{1}{2} \left( \{\beta \beta^*\}_\tau - \{ \alpha
\alpha^*\}_\tau \right)
              \end{pmatrix} \ .
\end{equation}
It is worth pointing out here that we have $\beta \beta^* = * \beta \wedge
\bar{\beta}$ which is true because $\Lambda^{0,2}$ is 1-dimensional. In other
words, the two diagonal entries only ``look" differently. 

With the above formulae (\ref{cliffordmultiplication}) we can now
write down the monopole equations (\ref{monopole equations}) with parameter
$\tau$ and as perturbation the imaginary-valued self-dual 2-form $\eta$ for the pair consisting of the spinor
$\Psi = (\alpha,\beta) \in \Gamma(X;\Lambda^{0,0}(E) \oplus \Lambda^{0,2}(E))$
and the connection ${A}$ in $E$:

\begin{equation}\label{un equations kaehler} 
\begin{split}
  	\bar{\partial}_{{A}} \alpha + \bar{\partial}^*_{{A}} \beta & = 0
\\
	F_{{A}}^{0,2} & =  \frac{1}{4} \{\beta \alpha^*\}_\tau + 4
	\eta^{0,2}
\\
  	- i \Lambda_g (F_{{A}}) & =  
	\frac{1}{8} \left\{ \alpha \alpha^* - *(\beta \wedge
	\bar{\beta})\right\}_\tau \	- i \Lambda_g(\eta)  
	\ .
\end{split}
\end{equation}
Indeed, the curvature equation of (\ref{monopole equations}) splits into four
equations according to the above splitting, but the two equations resulting
from the diagonal entries are equivalent, and, using that
$\overline{F_{{A}}^{0,2}} = -F_{{A}}^{2,0}$ (here again, $^-$ denotes
the complex-conjugation on the forms and the hermitian adjoint on $\End(E)$),
the two off-diagonal equations also prove to be equivalent.

\subsection{Decoupling phenomena, moduli spaces for $b_2^+(X) >  1$ and holomorphic 2-forms}
As mentioned before a lot of the analysis of the classical monopole equations on K\"ahler surfaces carries over to our situation. Before we consider the perturbed monopole equations we shall first draw some intermediate conclusions from the unperturbed monopole equations.
In particular there is a decoupling result completely analogous to the classical situation, interpreting monopoles as `vortices', c.f. also \cite{BG}, \cite{T}. 

\begin{prop}\label{decoupling}
Let $X$ be a K\"ahler surface. Suppose that the configuration $(\Psi,{A})
\in \Gamma(X;S^+_\mathfrak{c} \tensor E) \times \mathscr{A}(E)$ solves the {\em
unperturbed} $U(n)$ monopole equations with parameter $\tau \in [0,1]$. If we
write the spinor as $\Psi = (\alpha,\beta)$ according to the decomposition
$S^+_\mathfrak{c} \tensor E \cong \Lambda^{0,0}(E) \oplus \Lambda^{0,2}(E)$ then
one of the following two statements holds:
\begin{enumerate}
\item The second factor of the spinor vanishes identically, $\beta \equiv 0$.
Furthermore the pair  $(\alpha,{A})$ satisfies the following
`Vortex-type' equations
\begin{equation}\label{vortexalpha}
\begin{split}
  \bar{\partial}_{{A}} \alpha & = 0 \\
  F_{{A}}^{0,2} & = 0 \\
  i \Lambda_g (F_A) & = - \frac{1}{8} \{ \alpha \alpha^*\}_\tau \ .
\end{split}
\end{equation}
\item The first factor of the spinor vanishes identically, $\alpha \equiv 0$.
Furthermore the pair $(\beta,{A})$ satisfies the following equations
\begin{equation}\label{vortexbeta}
\begin{split}
  \bar{\partial}_{{A}}^* \beta & = 0 \\
  F_{{A}}^{0,2} & = 0 \\
  i \Lambda_g (F_A) & = + \frac{1}{8} \{\beta \beta^*\}_\tau \ .
\end{split}
\end{equation}
\end{enumerate}
\end{prop}
{\em Proof:} 
Using the first two of the monopole equations (\ref{un
equations kaehler}) we get:
\begin{equation*}
\dbar_A \dbar_A^* \beta  = - \dbar_A \dbar_A \alpha =
-F_A^{0,2} \alpha = -\frac{1}{4} \{\beta \alpha^*\}_\tau \alpha \ .
\end{equation*}
We take the inner product with $\beta$ to get now:
\begin{equation*}
\begin{split}
 \klammer{\beta,\dbar_A \dbar_A^* \beta} & = -\frac{1}{4}
\klammer{\beta, \geklatau{\beta \alpha^*} \alpha} \\
& = -\frac{1}{4} \klammer{\abs{\beta}^2 \abs{\alpha}^2 - \frac{1-\tau}{n}
\klammer{\beta, \tr(\beta \alpha^*) \alpha}} \\
& \leq \klammer{-\frac{1}{4} + \frac{1-\tau}{4n}} \abs{\alpha}^2 \abs{\beta}^2
\\
& \leq 0 \ .
\end{split}
\end{equation*}
Here we have used the Cauchy-Schwarz inequality, noting also that
$\abs{\tr(\beta \alpha^*)} \leq \abs{\beta \alpha^*} = \abs{\beta}
\abs{\alpha}$. Integrating now the latter inequality over the whole manifold $X$
yields the following:
\begin{equation}\label{conclusion inequality}
 0 \ \leq \norm{\dbar_A^* \beta}^2 \ \leq \klammer{-\frac{1}{4} +
	\frac{1-\tau}{4n}} \int_X \abs{\alpha}^2 \abs{\beta}^2 vol_g \leq 0
\end{equation}
Thus we get $\dbar_A^* \beta = 0$ and from the Dirac equation also
$\dbar_A \alpha = 0$. If further we have $\tau > 1-n$ then we see from the
last inequality that at any point of the manifold $X$ we have  $\alpha = 0$ or
$\beta = 0$. But we have $ 0 = \dbar_A^* \dbar_A \alpha =
\Delta_{\dbar_A} \alpha$, and because $\Delta_{\dbar_A}$ is an elliptic
second order operator with scalar symbol it follows from Aronaszajin's theorem
\cite{A} that solutions to $\Delta_{\dbar_A} \alpha = 0$ satisfy a unique
continuation theorem. Similarly we have $ 0 = \dbar_A \dbar_A^* \beta =
\Delta_{\dbar_A} \beta$, so the same holds for $\beta$. Therefore, if one
of $\alpha$ or $\beta$ vanishes on an open subset of $X$, then it vanishes on
the whole of $X$. The conclusions now follow from (\ref{un equations kaehler}).
\qed
\begin{remark}
If $\tau \neq 0$ the moduli space $M_{\mathfrak{c},E}(\tau,0)$
can only contain either solution with $\alpha \neq 0$ or with $\beta \neq 0$. 
This follows from taking the trace of the third equation of
(\ref{un equations kaehler}) and then integrating it over the whole manifold.
The left hand term yields then the topological quantity $- 2 \pi \, \langle
c_1(E) \smile \eckklammer{\omega_g},[X] \rangle$. 
\end{remark}

On a K\"ahler surface we have $\Delta = 2
\Delta_{\dbar}$, just reflecting again the compatibility between the complex structure and the Riemannian
metric. Therefore the harmonic differential forms are also $\dbar$-harmonic and vice versa. In
particular, we get the following decomposition from the Hodge-theorem:
\begin{equation}\label{h2_kaehler}
H^2_{dR}(X;\C) = H^{2,0}_{\bar{\partial}}(X) \oplus H^{1,1}_{\bar{\partial}}(X)
\oplus H^{0,2}_{\bar{\partial}}(X) \ .
\end{equation}

\begin{corollary}
If there are solutions to the unperturbed $U(n)$-monopole equations associated
to the data $(\mathfrak{c},E)$ and to the parameter $\tau \in [0,1]$, then the
image $c_1^\R(E)$ in real (complex) cohomology of the first Chern-class $c_1(E)
\in H^2(X;\Z)$ is of type $(1,1)$ according to the above decomposition
(\ref{h2_kaehler}).
\end{corollary}
{\em Proof:}
Under these conditions there is a connection ${A}$ on $E$ with
$F_A^{0,2} = 0 = F_A^{2,0}$. From the Chern-Weil formula we have that 
$\frac{-1}{2 \pi i} \eckklammer{\tr (F_A)} = c_1^\R(E)$. There is a 1-form
$\lambda$ such that $\omega:=\tr(F_A) - \dbar \lambda$ is $\dbar$-harmonic
and this class also represents $c_1^\R(E)$. We have $\omega^{2,0} = 0$ and
$\omega^{0,2} = \dbar \lambda$. But a class is $\dbar$ - harmonic if and only 
each component according to $\Omega^{p+q}(X) = \oplus \Omega^{p,q}(X)$ is
$\dbar$ - harmonic. But then the harmonic form $\omega^{2,0} = \dbar \lambda$
must be zero, as it is a $\dbar$ - exact form also. \qed

In the classical theory a common perturbation of the monopole equations was to
perturb with imaginary-valued self-dual 2-forms $\eta$ such that $\eta^{2,0}$ is
a holomorphic form \cite{W} \cite{Bq}. There are such forms with $\eta^{2,0}
\neq 0$ precisely if $b_2^+(X) > 1$. We will now consider this type of
perturbation in the general case of $U(n)$ monopoles even though these perturbations are not enough to get generic regularity of the moduli space in the case $n > 1$.  However, it will turn out that the moduli spaces perturbed in this way are empty in the case $n > 1$ as soon as the perturbing form $\eta$  is non-zero. 

If the unperturbed $U(n)$ monopole moduli space is empty then any invariant 
derived by the scheme `evaluation of cohomology classes on the fundamental cycle of the moduli space' should be zero. Indeed, that kind of invariant would be defined with a `generic' moduli space, i.e. one which is cut out transversally by the suitably perturbed monopole equations. An empty moduli space is always generic. Thus if there is a non-trivial invariant derived from some generic moduli space then the associated unperturbed moduli space may be not generic, but it could not be empty. Therefore it is natural to consider topological data $(\data)$ only for situations where the unperturbed $U(n)$ monopole moduli spaces are a priori non-empty. As we have seen, this can only be the case if the first Chern-class $c_1^\R(E)$ is of type $(1,1)$ according to the decomposition (\ref{h2_kaehler}). Therefore we shall include this hypothesis to the next two results, the following theorem and its corollary:

\begin{theorem}
Let $X$ be a K\"ahler surface and let $E$ be a bundle such that its first
Chern-class $c_1^\R(E)$ is of type $(1,1)$. Further let $\eta$ be an
imaginary-valued 2-form with $\eta^{2,0}$ holomorphic. Then the $U(n)-$
monopole equations (\ref{un equations kaehler}) associated to the data
$(\mathfrak{c},E)$, to the perturbation form $\eta$, and to the parameter $\tau \in (0,1]$ are equivalent to the
following system of equations:
\begin{equation}\label{equations kaehler split}
\begin{split}
  	\dbar_{{A}} \alpha & = 0 \\
	\dbar_{{A}}^* \beta & = 0 \\
	F_{{A}}^{0,2} & = 0 \\
	\frac{1}{4} \{\beta \alpha^*\}_\tau & = \eta^{0,2} \\
	-i \Lambda_g(F_{{A}}) & = \frac{1}{8} \{ \alpha \alpha^* - \beta
	\beta^* \}_\tau - i \Lambda_g(\eta) \\
\end{split}
\end{equation}
\end{theorem}
{\em Proof: } 
We will derive the following formula for a
solution $((\alpha,\beta),A))$ to the $U(n)-$monopole equations (\ref{un
equations kaehler}) with parameter $\tau$ and perturbation $\eta$:
\begin{equation}\label{formula}
\begin{split}
0 =  & 4 \, \norm{F_A^{0,2}}^2_{L^2(X)} \ 
+ 4 \, \frac{1-\tau}{\tau n} \, \norm{\tr F_A^{0,2}}^2_{L^2(X)}
+ \norm{\dbar_A^*
\beta}^2_{L^2(X)} \\
	& - \frac{4}{\tau} \, \langle 2 \pi i \ [\eta^{2,0}] \smile c_1(E) , [X]
	\rangle \
	 \ .
\end{split}
\end{equation}
The conclusion then clearly follows as the topological term vanishes by assumption.
\\

Provided that we have $\tau \neq 0$ the endomorphism $\beta \alpha^*$ can be
expressed as
\begin{equation}\label{expression of ab}
\begin{split}
\beta \alpha^* & = \geklatau{\beta \alpha^*} \ + \frac{1-\tau}{n} \tr (\beta
		\alpha^*) \\
		& = \geklatau{\beta \alpha^*} \ + \frac{1-\tau}{\tau n} \tr
		(\geklatau{\beta \alpha^*}) \\
		& = 4 \, F_A^{0,2} - 4\,  \eta^{0,2} + 4 \, \frac{1-\tau}{n \tau} \,
 			\tr(F_A^{0,2}) - 4 \, \frac{1-\tau}{\tau} \, \eta^{0,2} 	\ ,
\end{split}
\end{equation}
where the last equation used the second of the monopole equations (\ref{un
equations kaehler}) and the trace of it. 

Again we get from the Dirac-equation that $
\dbar_A \dbar_A^* \beta + F_A^{0,2} \alpha = 0 $, 
so that after taking the pointwise inner-product with $\beta$ and using the
above equation (\ref{expression of ab}) we get: 
\begin{equation}\label{to integrate}
\begin{split}
0 &  = \klammer{\beta, F_A^{0,2} \alpha} \ + \klammer{\beta, \dbar_A
		\dbar_A^* \beta } \\
	& = \klammer{ \beta \alpha^* , F_A^{0,2} } \ + \klammer{\beta, \dbar_A \dbar_A^* \beta } \\ 
	& = 4 \, \abs{F_A^{0,2}}^2 \ - 4 \, \klammer{\eta^{0,2}, F_A^{0,2} } \\
	& \ \ \ \ \ 	+ 4 \frac{1-\tau}{n \tau} \abs{ \tr(F_A^{0,2})}^2 \ - 4  \, \frac{1-\tau}{\tau} \, 
		\klammer{\eta^{0,2}, F_A^{0,2}} + \,  \klammer{\beta, \dbar_A
		\dbar_A^* \beta } \\
	& = 4  \, \abs{F_A^{0,2}}^2 \  + 4 \, \frac{1-\tau}{n \tau} \abs{ \tr(F_A^{0,2})}^2  -
\frac{4}{\tau} \klammer{\eta^{0,2}, F_A^{0,2} } + \,  \klammer{\beta, \dbar_A
		\dbar_A^* \beta }
\end{split}
\end{equation}
As the next step we will integrate this whole equation over $X$. Beforehand we shall remark that  $\eta^{2,0}$ is closed, and therefore the following integral is of topological nature:
\begin{equation}\label{topological term}
\begin{split}
\int_X \klammer{\eta^{0,2}, F_A^{0,2} } \, vol_g & = \int_X \bar{\eta^{0,2}} \wedge * \tr(F_A^{0,2}) \\
	& = - \, \int_X \eta^{2,0} \wedge \tr(F_A^{0,2}) \\
	& = \, 2 \pi i \ \langle  [\eta^{2,0}] \smile c_1(E) , [X]
	\rangle
\end{split}
\end{equation}
With this said the integral of the formula (\ref{to integrate}) clearly yields the above formula (\ref{formula}).

\qed

\begin{corollary}\label{final conclusion}
Let $X$ be a K\"ahler surface with $b_2^+(X) > 1$ and let $E$ be a bundle such
that its first Chern-class $c_1^\R(E)$ is of type $(1,1)$. Then for any
self-dual imaginary valued 2-form $\eta$ with $\eta^{2,0}$ holomorphic and {\em
non-zero} and constant $\tau \in (0,1]$ the moduli space $M_{\mathfrak{c},E}(\eta,\tau)$ is empty. 
\end{corollary}
{\em Proof:}
  Under the given hypothesis the preceeding theorem implies that 
\begin{equation}
\label{impossible}
\geklatau{ \beta \alpha^*} = 4 \eta^{0,2} \ \id_E \ .
\end{equation}
But using the definition of $\geklatau{\beta \alpha^*}$ it is a pure matter
of linear algebra to check that for $\eta^{0,2} \neq 0$ this is impossible if $n
\geq 2$, because the left hand side of the equation (\ref{impossible}) can
never be a mutliple of the identity, unless $\alpha = 0$ or $\beta = 0$.
\qed


\begin{thebibliography}{99999}
\bibitem{A} N. Aronaszajin, {\em A unique continuation theorem for solutions
to elliptic partial differential equations or inequalities of the second
order}, {Journal de Math\'ematiques Pures et Appliqu\'ees} (9), 36, (1957),
235-249.
\bibitem{BG} S. Bradlow, O. Garcia-Prada, {\em Non-abelian monopoles and
vortices}, Lecture Notes in Pure and Appl. Math. 184 (1997), 567-589.
\bibitem{Bq} O. Biquard,  {\em Les \'equations de Seiberg-Witten sur  une
surface
complexe non-k\"ahlerienne},  Comm. Anal. Geom. 6, No.1 (1998), 173-197.
\bibitem{D} S. Donaldson, {\em Polynomial invariants for smooth
four-manifolds}, Topology 2, No. 3, (1990), 257-315.
\bibitem{DK} S. Donaldson, P. Kronheimer, {\em The Geometry of
Four-Manifolds}, {Oxford Mathematical Monographs} (1990).
\bibitem{FL1} P. Feehan, T. Leness, {\em A general SO(3)-monopole
cobordism formula relating Donaldson and Seiberg-Witten invariants},
{Preprint (2002)}, arxiv:math.DG/0203047.
\bibitem{FL2} P. Feehan, T. Leness, {\em On Donaldson and Seiberg-Witten
invariants}, {Proc. Sympos. Pure Math. 71} (2001) 237-248.
\bibitem{FL3} P. Feehan, T. Leness, {\em SO(3) Monopoles, Level-One
Seiberg-Witten Moduli Spaces, and Witten's Conjecture in Low Degrees},
{Topology and its Applications}, to appear.
\bibitem{FL4} P. Feehan, T. Leness, {\em PU(2)-monopoles. I: regularity,
Uhlenbeck compactness, and transversality}, {Journal of Differential Geometry},
49 (1998), 265-410. 
\bibitem{FL5} P. Feehan, T. Leness, {\em Witten's conjecture for four-manifolds of simple type}, {arxiv:math.DG/0609530} (2006). 
\bibitem{GP} O. Garcia-Prada, {\em A direct existence proof for the vortex
equations over a compact Riemann surface}, Bull. London Math. Soc. 26 (1994),
88-96. 
\bibitem{Hi} N. Hitchin, {\em Harmonic spinors}, {Advances in Mathematics
 14}, (1974) 1-55.
\bibitem{K} P. Kronheimer, {\em Four-manifold invariants from higher rank
bundles}, {Journal of Differential Geometry}, No. 70 (2005), 59-112.
\bibitem{KM} P. Kronheimer, T. Mrowka, {\em The genus of embedded surfaces
in the projective plane}, {Math Research Letters 1} , No. 1, (1994) 796-808.
\bibitem{KM2} P. Kronheimer, T. Mrowka, {\em Embedded surfaces and the
structure of Donaldson's polynomial invariants}, Journal of Differential
Geometry 41, No.3, (1995), 572-734.
\bibitem{MM} M. Marino, G. Moore, {\em The Donaldson-Witten function for
gauge groups of rank larger than one}, Commun. Math. Phys. 199 (1998), 25-69.
\bibitem{M} J. Morgan, {\em The Seiberg-Witten Equations and Applications to
the Topology of smooth Four-Manifolds}, {Math. Notes, Princeton Univ. Press},
1996
\bibitem{MS} J. Milnor, J. Stasheff, {\em Characteristic Classes}, {Annals
of Mathematics Studies, Princeton University Press}, (1974).
\bibitem{N} S. Nicolaescu, {\em Notes on Seiberg-Witten Theory}, {GSM Vol.
28, American Mathematical Society}, (2000).
\bibitem{OT} C. Okonek, A. Teleman, {\em Master Spaces and the Coupling
Principle: From Geometric Invariant Theory to Gauge Theory}, Commun. Math. Phys.
205 (1999) 437-458.
\bibitem{OT2} C. Okonek, A. Teleman, {\em Quaternionic Monopoles}, {Comm.
Math. Phys.}, Vol. 180, No. 2 (1996), 363-388.
\bibitem{OT3} C. Okonek, A. Teleman, {\em The coupled Seiberg-Witten
equations, vortices and moduli spaces of stable pairs}, Int. J. Math., Vol. 6,
No. 6 (1995), 893-910.
\bibitem{OT4} C. Okonek, A. Teleman, {\em Recent developments in
Seiberg-Witten theory and complex geometry}, {Several complex variables
(Berkeley, 1995-1996)}, Math. Sci. Res. Inst. Publ. 37, Cambridge University
Press, 1999, 391-428.
\bibitem{PT} V. Pidstrigach, A. Tyurin, {Localization of Donaldson
invariants along the Seiberg-Witten classes}, preprint dg-ga/9507004
\bibitem{S} N. Steenrod, {\em The Topology of fibre bundles}, {Princeton University Press} (1951).
\bibitem{T} A. Teleman, {\em Non-abelian Seiberg-Witten theory and stable
oriented pairs}, {Int. J. of Math.} No.4 (1997) 507-535.
\bibitem{T2} A. Teleman, {\em The moduli space of PU(2)-monopoles}, {Asian
J. of Math.} 
\bibitem{T3} A. Teleman, {\em Introduction \`a la Th\'eorie de
Jauge}, {Lecture Notes}, available at \begin{verbatim} http://www.cmi.univ-mrs.fr/~teleman
\end{verbatim}
\bibitem{W} E. Witten, {\em Monopoles and Four-manifolds}, {Math Research
Letters 1} No.1 (1994) 809-822.
\bibitem{We} R. Wells, {\em Differential Analysis on Complex Manifolds},
{GTM No. 65, Springer}, (1980). 
\bibitem{Z1} R. Zentner, {\em PhD thesis}, available online at \begin{verbatim} http://www.cmi.univ-mrs.fr/~zentner/thesis.pdf \end{verbatim}
\bibitem{Z} R. Zentner, {A vanishing result for a Casson-type instanton invariant}, preprint (2009), arXiv:0911.2772
\end{thebibliography}
\end{document}